\documentclass{article}
\usepackage{amssymb}
\usepackage{amsfonts}
\usepackage{amsmath,amsthm,amssymb,latexsym,amscd}
\usepackage{epsf}
\usepackage[pdftex]{graphicx}
\usepackage{tikz}
\newtheorem{lemma}{Lemma}
\newtheorem{theorem}{Theorem}
\newtheorem{proposition}{Proposition}
\newtheorem{remark}{Remark}
\newtheorem{definition}{Definition}
\newtheorem{hypothesis}{Hypothesis}
\title{Cohesive dynamics and brittle fracture}
\author{Robert Lipton\thanks{Department of Mathematics,
        Louisiana State University,
        Baton Rouge, LA 70803,
        ({\tt lipton@math.lsu.edu}).}}
\begin{document}
\maketitle
\begin{abstract} 
We formulate a nonlocal cohesive model for calculating the deformation state inside a cracking body. In this model a more complete set of physical properties including elastic and softening behavior are assigned to each point in the medium. We work within the small deformation setting and use the peridynamic formulation. Here strains are calculated as difference quotients. The constitutive relation is given by a nonlocal cohesive law relating force to strain. At each instant of the evolution we identify a process zone where strains lie above a threshold value.  Perturbation  analysis shows that jump discontinuities within the process zone can become unstable and grow.  We derive an explicit  inequality  that shows that the size of the process zone is controlled by the ratio given by the length scale of nonlocal interaction divided by the characteristic dimension of the sample. The process zone is shown to concentrate on a set of zero volume in the  limit where the length scale of nonlocal interaction vanishes with respect to the size of the domain. In this limit the dynamic evolution is seen to have bounded linear elastic energy and Griffith surface energy. The limit dynamics corresponds to the simultaneous evolution of linear elastic displacement and the fracture set across which the displacement is discontinuous. We conclude illustrating how the approach developed here can be applied to limits of dynamics associated with other energies that $\Gamma$- converge to the Griffith fracture energy.
\end{abstract}
\begin{flushleft} 
\end{flushleft}

\begin{flushleft}
\end{flushleft}
%
%
%\begin{flushleft}{\bf AMS Subject Classifications:} 
%\end{flushleft}
%
\pagestyle{myheadings}
\markboth{R. LIPTON}{Cohesive Dynamics and Brittle Fracture}
\setcounter{equation}{0} \setcounter{theorem}{0} \setcounter{lemma}{0}\setcounter{proposition}{0}\setcounter{remark}{0}\setcounter{definition}{0}\setcounter{hypothesis}{0}

\section{Introduction}
\label{introduction}

Dynamic brittle fracture is a multiscale phenomenon operating across a wide range of length and time scales. Contemporary approaches to brittle fracture modeling can be broadly characterized as bottom-up and top-down. Bottom-up approaches take into account the discreteness of fracture  at the smallest length scales and are expressed through lattice models. This approach has provided insight into the dynamics of the fracture process \cite{6,16,17,21}. Complementary to  the bottom-up approaches are top-down computational approaches using cohesive surface elements \cite{Cox}, \cite{14}, \cite{22}, \cite{Remmers}. In this formulation the details of the process zone are collapsed onto an interfacial element with a force traction law given by the cohesive zone model  \cite{Barenblatt}, \cite{Dougdale}. Cohesive surfaces have been applied within the extended finite element method \cite{5}, \cite{Duarte}, \cite{18} to minimize the effects of mesh dependence on free crack paths. Higher order multi-scale cohesive surface models involving excess properties and differential momentum balance are developed in \cite{WaltonSendova}. Comparisons between different cohesive surface models are given in \cite{Falk}. More recently variational approaches to brittle fracture based on quasi-static evolutions of global minimizers of Grifffith's fracture energy have been developed  \cite{FrancfortMarigo}, \cite{BourdinFrancfortMarigo}, and \cite{FrancfortLarsen}. Phase field approaches have also been developed to  model brittle fracture evolution from a continuum perspective \cite{BourdinFrancfortMarigo}, \cite{BourdinLarsenRichardson}, \cite{Miehe}, \cite{Hughes}, \cite{Ortiz}, \cite{Wheeler}. In the phase field approach a second field is introduced to interpolate between cracked and undamaged elastic material. The evolution of the phase field is used to capture the trajectory of the crack.  A concurrent development is the emergence of the peridynamic formulation introduced in \cite{Silling1} and \cite{States}. Peridynamics is a  nonlocal formulation of continuum mechanics expressed in terms of displacement differences as opposed to spatial derivatives of the displacement field.  These features provide the ability to simultaneously simulate kinematics involving both smooth displacements and defect evolution.  Numerical simulations based on peridynamic modeling exhibit the formation and evolution of sharp interfaces associated with defects and  fracture \cite{Bobaru1},  \cite{BhattacharyaDyal}, \cite{Foster}, \cite{Bobaru2}, \cite{SillingBobaru}, \cite{SillingAscari2}, \cite{WecknerAbeyaratne}.  In an independent development nonlocal formulations have  been introduced for modeling the passage from discrete to continuum limits of energies for quasistatic fracture models \cite{Alicandro},  \cite{Buttazzo}, \cite{Braides}, \cite{BraidesGelli},  for smeared crack models \cite{Negri} and for image processing \cite{Gobbino1} and \cite{Gobbino3}.  A complete review of contemporary methods is beyond the scope of this paper however the reader is referred to \cite{Agwai}, \cite{Bazant},  \cite{BelitchReview}, \cite{Bouchbinder}, \cite{BourdinFrancfortMarigo}, \cite{Braides1}, \cite{Braides2} for a more complete guide to the literature.

In this paper we formulate a nonlocal, multi-scale, cohesive continuum model for assessing the deformation state inside a cracking body. This model is expressed using the peridynamic formulation introduced in \cite{Silling1}, \cite{States}. Here strains are calculated as difference quotients of displacements between two points $x$ and $y$. In this approach the force between two points $x$ and $y$ is related to the strain through a nonlinear cohesive law that  depends upon the magnitude and direction of the strain. The forces are initially elastic for small strains and soften beyond a critical strain.  We introduce the dimensionless length scale  $\epsilon$  given by the ratio of the length scale of nonlocal interaction to the characteristic length of the material sample $D$.  Working in the new rescaled coordinates the nonlocal interactions between $x$ and its neighbors $y$  occur within a horizon of radius $\epsilon$ about $x$ and the characteristic length of $D$ is taken to be unity. This neighborhood of $x$ is denoted by $\mathcal{H}_\epsilon(x)$. 

To define the potential energy we first assume the deformation $z$ is given by $z(x)=u(x)+x$ where $u$ is the displacement field. The strain between two points $x$ and $y$ inside $D$ is given by
\begin{eqnarray}
\mathcal{S}=\frac{|z(y)-z(x)|-|y-x|}{|y-x|}.
\label{bigdefstrain}
\end{eqnarray}
In this treatment we assume small deformation kinematics and the displacements are small (infinitesimal) relative to the size of the body $D$. Under this hypothesis \eqref{bigdefstrain} is linearized and the strain is given by  $$\mathcal{S}=\frac{u(y)-u(x)}{|y-x|}\cdot e,$$ where  $e=\frac{y-x}{|y-x|}$. Both two and three dimensional problems will be considered and the dimension is denoted by $d=2,3$.
The cohesive model is characterized through a nonlocal potential $$W^\epsilon(\mathcal{S},y-x),$$ associated with points $x$ and $y$.  The associated energy density is obtained on integrating over $y$ for $x$ fixed and is given by
\begin{eqnarray}
{\bf{W}}^\epsilon(\mathcal{S},x)=\frac{1}{V_d}\int_{\mathcal{H}_\epsilon(x)}\,W^\epsilon(\mathcal{S},y-x)\,dy
\label{densityy}
\end{eqnarray}
where $V_d=\epsilon^d\omega_d$ and $\omega_d$ is the (area) volume of the unit ball in dimensions $d=(2)3$. The potential energy of the body  is given by
\begin{eqnarray}
PD^\epsilon(u)=\int_{D}\,{\bf{W}}^\epsilon(\mathcal{S},x)\,dx
\label{the peridynamicenergy}
\end{eqnarray}

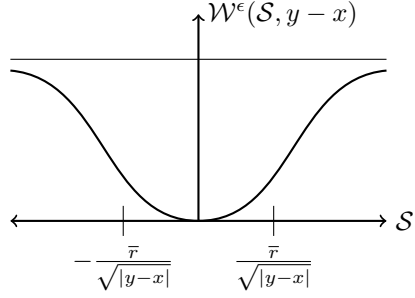
\begin{figure} 
\centering
\begin{tikzpicture}[xscale=1,yscale=1]
\draw [<->,thick] (0,2.75) -- (0,0) -- (2.5,0);
\draw [->,thick] (0,0) -- (-2.5,0);
\draw [-] (-2.5,2.15) -- (2.5,2.15);
\draw [-,thick] (0,0) to [out=0,in=-175] (2.5,2);
\draw [-,thick] (-2.5,2) to [out=-5,in=180] (0,0);
\draw (1,-0.2) -- (1, 0.2);
\draw (-1,-0.2) -- (-1, 0.2);
\node [below] at (1,-0.2) {$\frac{\overline{r}}{\sqrt{|y-x|}}$};
\node [below] at (-1,-0.2) {$-\frac{\overline{r}}{\sqrt{|y-x|}}$};
\node [right] at (2.5,0) {$\mathcal{S}$};
\node [right] at (0,2.75) {$\mathcal{W}^\epsilon(\mathcal{S},y-x)$};
\end{tikzpicture} 
\caption{\bf Cohesive potential as a function of $\mathcal{S}$ for $x$ and $y$ fixed.}
\label{ConvexConcave}
\end{figure}

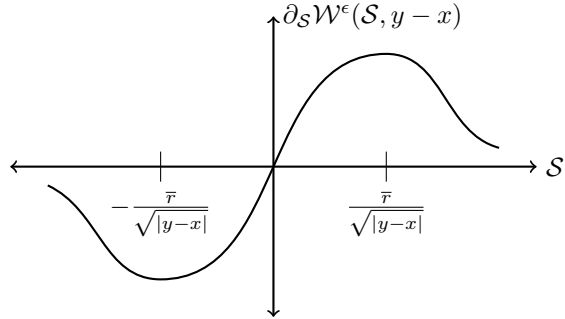
\begin{figure} 
\centering
\begin{tikzpicture}[xscale=1,yscale=1]
\draw [<->,thick] (0,2) -- (0,-2);
\draw [<->,thick] (-3.5,0) -- (3.5,0);
\draw [-,thick] (-3,-0.25) to [out=-25,in=180] (-1.5,-1.5) to [out=0,in=180] (1.5,1.5)
           to [out=0,in=165] (3,0.25);
\draw (1.5,-0.2) -- (1.5, 0.2);
\draw (-1.5,-0.2) -- (-1.5, 0.2);
\node [below] at (1.5,-0.2) {$\frac{\overline{r}}{\sqrt{|y-x|}}$};
\node [below] at (-1.5,-0.2) {$-\frac{\overline{r}}{\sqrt{|y-x|}}$};
\node [right] at (3.5,0) {$\mathcal{S}$};
\node [right] at (0,2.0) {$\partial_{\mathcal{S}}\mathcal{W}^\epsilon(\mathcal{S},y-x)$};
\end{tikzpicture} 
\caption{{\bf Cohesive relation between force and  strain for $x$ and $y$ fixed.}}
 \label{SofteningBond}
\end{figure}

We  introduce the class of potentials associated with a cohesive force that is initially elastic and then softens after a critical strain. These potentials are of the generic form given by
\begin{eqnarray}
W^\epsilon(S,y-x)=|y-x|\mathcal{W}^\epsilon(\mathcal{S},y-x),
\label{potentialdensity1a}
\end{eqnarray}
where $\mathcal{W}^\epsilon(\mathcal{S},y-x)$ is the peridynamic potential per unit length associated with $x$ and $y$ given by
\begin{eqnarray}
\mathcal{W}^\epsilon(\mathcal{S},y-x)=\frac{1}{\epsilon}J^\epsilon\left(|y-x|\right)\left(\frac{1}{|y-x|}f\left(|y-x|\mathcal{S}^2\right)\right).
\label{potentialdensity2a}
\end{eqnarray}
These potentials are of a general form and are associated with potential functions $f:[0,\infty)\rightarrow\mathbb{R}$  that are positive, smooth and concave with the properties
\begin{eqnarray}
\lim_{r\rightarrow 0^+}\frac{f(r)}{r}=f'(0)>0,&&\lim_{r\rightarrow\infty}f(r)=f_\infty <\infty.
\label{properties}
\end{eqnarray}
The composition of $f$ with $|y-x|\mathcal{S}^2$ given by \eqref{potentialdensity2a} delivers the convex-concave dependence of $\mathcal{W}^\epsilon(\mathcal{S},y-x)$ on $\mathcal{S}$ for fixed values of $x$ and $y$, see Figure \ref{ConvexConcave}.  Here $J^\epsilon(|y-x|)$ is used to prescribe the influence of separation length $|y-x |$ on the force between $x$ and $y$ with $0\leq J^\epsilon(|y-x|)<M $ for $0\leq |y-x|< \epsilon$ and $J^\epsilon(|y-x|)=0$ for $\epsilon\leq |y-x |$. For fixed $x$ and $y$ the inflection point for the potential energy \eqref{potentialdensity2a} with respect to the strain $\mathcal{S}$   is given by $\overline{r}/\sqrt{|y-x|}$, where $\overline{r}$ is the inflection point for the function $r:\rightarrow f(r^2)$,  see Figure \ref{ConvexConcave}. This choice of potential delivers an initially elastic and then softening constitutive law for the force per unit length along the direction $e$ given by
\begin{eqnarray}
\hbox{\rm force per unit length}=\partial_\mathcal{S} \mathcal{W}^\epsilon(\mathcal{S},y-x)=\frac{2}{\epsilon}\left(J^\epsilon(|y-x|)f'\left(|y-x|\mathcal{S}^2\right)\mathcal{S}\right).
\label{forcestate}
\end{eqnarray}
The force between points $y$ and $x$ begins to drop when the strain $\mathcal{S}$ exceeds the critical  strain 
\begin{eqnarray}
|\mathcal{S}|>\frac{\overline{r}}{\sqrt{|y-x|}}=\mathcal{S}_c,
\label{sqrtsingularity}
\end{eqnarray}
see Figure \ref{SofteningBond}. This is the same singularity strength  associated with a strain concentration in the vicinity of a crack tip as in the classic theory of brittle fracture  \cite{Freund}.    A future goal will be to inform the cohesive model introduced here by resolving atomistic or molecular dynamics across smaller length scales.

We apply the principle of least action to recover the cohesive equation of motion describing the state of displacement inside the body $D\subset\mathbb{R}^d$  given by
\begin{eqnarray}
\rho\partial_{tt}^2 u(t,x)=2\frac{1}{V_d}\int_{\mathcal{H}_\epsilon(x)}\,\partial_\mathcal{S} \mathcal{W}^\epsilon(\mathcal{S},y-x)\frac{y-x}{|y-x|}\,dy+b(t,x)
\label{eqofmotion}
\end{eqnarray}
where $\rho$ is the density and $b(t,x)$ is the body force. This is a well posed formulation in that existence and uniqueness (within a suitable class of evolutions) can be shown, see Section \ref{sec2} and Theorem \ref{existenceuniqueness} of Section \ref{EE}.
In this model a more complete set of physical properties including elastic and softening behavior are assigned to each point in the medium. Here each point in the domain is connected to its neighbors by a cohesive law  see Figure \ref{SofteningBond}.   We define the {\em process zone} to be the collection of points $x$ inside the body $D$ associated with peridynamic neighborhoods $\mathcal{H}_\epsilon(x)$ for which the strain $\mathcal{S}$ between $x$ and $y$  exceeds a threshold value for a sufficiently large proportion of points $y$ inside $\mathcal{H}_\epsilon(x)$. Here the force vs. strain law departs from linear behavior when the strain exceeds the threshold value. The mathematically precise definition of the process zone is given in Section \ref{sec4}, see Definition  \ref{processZone}. In this model the {\em fracture set} is  associated with peridynamic neighborhoods $\mathcal{H}_\epsilon(x)$ with strains  $|\mathcal{S}|>\mathcal{S}_c$ for which the force vs. strain law begins to soften and is defined in Section \ref{sec4}, see Definition \ref{Fractureset}.
The nonlinear elastic--softening behavior put forth in this paper is similar to the ones used in cohesive zone models \cite{Dougdale}, \cite{Barenblatt}. However for this model the dynamics selects whether a material point lies inside or outside the process zone.  The principle feature of the cohesive dynamics model introduced here is that the evolution of the process zone together with the the fracture set is goverened  by one equation   consistent with Newton's second law given by \eqref{eqofmotion}. This is a characteristic feature of peridynamic models  \cite{Silling1}, \cite{States} and  lattice formulations for fracture evolution \cite{6,16,17,21}.

In this paper the goal is to characterize the size of the process zone  for  cohesive dynamics as a function of domain size and the length scale of the nonlocal forces. The second focus is to identify properties of the distinguished limit evolution for these models in the limit of vanishing non-locality as characterized by the $\epsilon\rightarrow 0$ limit. For the model introduced here the  parameter that controls the size of the process zone is given by the radius of the horizon $\epsilon$. We derive an explicit  inequality  that shows that the size of the process zone is controlled by the horizon radius. Perturbation analysis shows that jump discontinuities within the process zone can become unstable and grow.
This analysis shows that {\em the horizon size $\epsilon$ for cohesive dynamics is a modeling parameter} that can be calibrated according to the size of the process zone obtained from experimental measurements.

 Further calculation shows that the volume of the process zone goes to zero with $\epsilon$ in the limit of vanishing non-locality, $\epsilon\rightarrow 0$.   Distinguished $\epsilon\rightarrow 0$  limits of cohesive evolutions are identified and are found to have both bounded linear elastic energy and Griffith surface energy.  Here the limit dynamics corresponds to the simultaneous evolution of linear elastic displacement and a fracture set across which the displacement is discontinuous.  Under suitable hypotheses it is found that for points in space-time away from the fracture set that the displacement field evolves according to the linear elastic wave equation. Here the linear wave equation provides a dynamic coupling between elastic waves and the evolving fracture path inside the media. The elastic moduli, wave speed and energy release rate for the evolution are explicitly recovered from moments of the peridynamic potential energy. These features are consistent with the asymptotic behavior seen in the convergence of solutions of the Barenblatt model to the Griffith model when cohesive forces confined to a surface act over a sufficiently short range \cite{MarigoTruskinovsky},
\cite{Willis}.

Earlier work has shown that linear peridynamic formulations recover the classic linear elastic wave equation in the limit of vanishing non-locality see \cite{EmmrichWeckner},  \cite{SillingLehoucq}. The convergence of linear peridynamics to Navier's equation in the sense of solution operators is demonstrated in \cite{MengeshaDu}.  Recent work shows that analogous results can be found for dynamic problems and fully nonlinear peridynamics \cite{LiptonJElast2014} in the context of antiplane shear. There distinguished $\epsilon\rightarrow 0$  limits of cohesive evolutions are identified and are found to have both bounded linear elastic energy and Griffith surface energy. It is shown that the limiting displacement evolves according to the linear elastic wave equation away from the crack set, see  \cite{LiptonJElast2014}. For large deformations, the connection between hyperelastic energies and the small horizon limits of nonlinear peridynamic energies is recently established in \cite{Bellido}. In the current paper both two and three dimensional problems involving multi-mode fracture are addressed. For these problems new methods are required to identify the existence of a limit dynamics as the length scale of nonlocal interaction $\epsilon$ goes to zero.  A crucial step is to establish a suitable notion of compactness for sequences of cohesive evolutions. The approach taken here employs nonlocal Korn inequalities  introduced in \cite{DuGunzbergerlehoucqmengesha}. This method is presented in Section \ref{CC}.  We conclude noting that the cohesive dynamics model introduced here does not have an irreversibility constraint and that the constitutive law  \eqref{forcestate}  applies at all times in the fracture evolution. However with this caveat in mind,  the nonlocal cohesive model offers new computational and analytical opportunities for understanding the effect of the process zone on fracture patterns.

In the next section  we write down the Lagrangian formulation for the cohesive dynamics and apply the principle of least action to recover the equation of motion. In that section it is shown that the nonlinear-nonlocal  cohesive evolution is a well posed initial boundary value problem. It is also shown that energy balance is satisfied by the cohesive dynamics. A formal stability analysis is carried out in Section \ref{sec3} showing that  jump discontinuities within the process zone can become unstable and grow, see Proposition \ref{Nuccriteria}.  In Section \ref{sec4} we provide a mathematically  rigorous inequality explicitly showing how the volume of the process zone for the cohesive evolutions is controlled by the length scale of nonlocal interaction $\epsilon$, see Theorem \ref{epsiloncontropprocesszone}. It is shown that the process zone concentrates on a set of zero volume in the limit, $\epsilon\rightarrow 0$, see Theorem \ref{bondunstable}. In Sections \ref{sec5} and \ref{sec6} we introduce suitable technical hypothesis  and identify the distinguished limit of the cohesive evolutions as $\epsilon\rightarrow 0$, see Theorem \ref{LimitflowThm}. It is shown that the dynamics can be expressed in terms of displacements  that satisfy the linear elastic wave equation away from the crack set, see Theorem \ref{waveequation}. These displacements are  shown to have  bounded bulk elastic and surface energy in the sense of Linear Elastic Fracture Mechanics (LEFM), see Theorem \ref{LEFMMThm}.   In Section \ref{EightIntro} we provide the mathematical underpinnings and proofs of the theorems. In Section \ref{Ninth}   we apply the approach developed here to examine limits of dynamics associated with other energies that $\Gamma$- converge to the Griffith fracture energy.  As an illustrative example we examine the Ambrosio-Tortorelli \cite{AT} approximation as applied  to the dynamic problem in \cite{BourdinLarsenRichardson} and \cite{LarsenOrtnerSuli}. 

\setcounter{equation}{0} \setcounter{theorem}{0} \setcounter{lemma}{0}\setcounter{proposition}{0}\setcounter{remark}{0}\setcounter{remark}{0}\setcounter{definition}{0}\setcounter{hypothesis}{0}

\section{Cohesive dynamics}
\label{sec2}

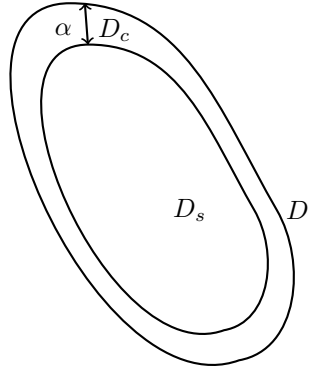
\begin{figure} 
\centering
\begin{tikzpicture}[xscale=1,yscale=1]
%\draw [<->,thick] (0,2) -- (0,-2);
\draw [<->,thick] (-1.01,2.2) -- (-1.051,2.75);
\draw [-,thick] (-1.25,2.75) to [out=0,in=120] (1.5,0) to [out=300,in=10] (1.0,-2.0) to  [out=200,in=180](-1.25, 2.75);
\draw [-,thick] (-1.0,2.2) to [out=0,in=120] (1.2,0) to [out=300,in=10] (0.8,-1.6) to  [out=200,in=180](-1.0, 2.2);
\node [right] at (1.5,0) {$D$};
\node [right] at (0,0) {$D_s$};
\node [right] at (-1.0,2.4) {$D_c$};
\node [left] at (-1.10,2.4) {$\alpha$};
%\node [right] at (0,2.0) {force};
\end{tikzpicture} 
\caption{{\bf Domain $D=D_c\cup D_s$.}}
 \label{Domains}
\end{figure}

We formulate the initial boundary value problem for the cohesive evolution. Since the problem is nonlocal the domain $D$ is split into a boundary layer called the constraint set $D_c$, and the interior $D_s$. To fix ideas the thickness of the boundary layer is denoted by $\alpha$ and  $2\epsilon<\alpha$ where $2\epsilon$ is the diameter of nonlocal interaction see Figure \ref{Domains}.  The boundary condition for the displacement $u$ is given by $u(t,x)=0$ for $x$ in $D_c$.  To  incorporate nonlocal boundary conditions we introduce the space $L^2_0(D;\mathbb{R}^d)$,  of displacements that are square integrable over $D$ and zero on $D_c$. 
The initial conditions for the cohesive  dynamics belong to $L^2_0(D;\mathbb{R}^d)$ and are  given by
\begin{eqnarray}
u(x,0)=u_0(x),\hbox{ and }u_t(x,0)=v_0(x).
\label{initialconditions}
\end{eqnarray}

We will investigate the evolution of the deforming domain for general initial conditions. These can include an initially un-cracked body or one with a preexisting system of cracks. For two dimensional problems the cracks are given by a system of curves of finite total length, while for three dimensional problems the crack set is given by a system of surfaces of finite total surface area. Depending on the dimension of the problem the displacement suffers a finite jump discontinuity across each curve or surface. The initial condition is specified by a crack set $K$ and displacement $u_0$. The strain  $\mathcal{E}u_0=(\nabla u_0+\nabla u_0^T)/2$ is defined off the crack set and the displacement $u_0$ can suffer jumps across $K$.  Griffith's theory of fracture
asserts that the energy necessary to produce a crack $K$ is proportional to the crack length (or surface area). For Linear Elastic Fracture Mechanics (LEFM) the total energy associated with bulk elastic and surface  energy is given by
\begin{eqnarray}
LEFM(u_0)=\int_{D}\left(2\mu|\mathcal{E}u_0|^2+\lambda|{\rm div}\,u_0|^2\right)\,dx+\mathcal{G}_c |K|,
\label{Gcrackenergy}
\end{eqnarray}
where $\mu$, $\lambda$ are the  the shear and Lam\'e  moduli and $\mathcal{G}_c$ is the critical energy release rate for the material. Here $|K|$ denotes the length or surface area of the crack.
In what follows we will assume that the  bulk elastic energy and surface energy of the initial displacement are bounded as well as the the initial velocity and displacement.
For future reference we describe initial data $u_0$ and $v_0$ that satisfy these conditions as {\em LEFM initial data} and we have the inequality between the peridynamic energy and the energy of Linear Elastic Fracture Mechanics given by
\begin{eqnarray}
PD^\epsilon(u_0)\leq LEFM(u_0),
\label{basicinequality}
\end{eqnarray}
when $\mu$, $\lambda$, and $\mathcal{G}_c$ are related to the nonlocal potentials according to  \eqref{calibrate1} and \eqref{calibrate2}.
This inequality is established in Section \ref{CC}, see \eqref{upperboundperi}.

In what follows we write $u(t,x)$ as $u(t)$ to expedite the presentation.
The cohesive dynamics is described by the Lagrangian
\begin{eqnarray}
L^\epsilon(u(t),\partial_t u(t),t)=K(\partial_t u(t))-PD^\epsilon(u(t))+U(u(t)),
\label{Lagrangian}
\end{eqnarray}
with 
\begin{eqnarray}
K(\partial_t u(t))&=&\frac{1}{2}\int_{D}\rho|\partial_t u(t,x)|^2\,dx, \hbox{ and }\nonumber\\
U(u(t))&=&\int_{D}b(t,x) u(t,x)\,dx,
\label{Components}
\end{eqnarray}
where $\rho$ is the mass density of the material and $b(t,x)$ is the body force density. 
The  initial conditions $u^\epsilon(0,x)=u_0(x)$ and $u^\epsilon_t(0,x)=v_0(x)$ are prescribed and the action integral for the peridynamic evolution is
\begin{eqnarray}
I^\epsilon(u)=\int_0^TL^\epsilon(u(t),\partial_t u(t),t)\,dt.
\label{Action}
\end{eqnarray}
The Euler Lagrange Equation for this system delivers the the  cohesive dynamics described by
\begin{eqnarray}
\rho u^\epsilon_{tt}&=&-\nabla PD^\epsilon(u^\epsilon)+b
\label{stationary}
\end{eqnarray}
where
\begin{eqnarray}
\nabla PD^\epsilon(u^\epsilon)=-\frac{2}{V_d}\int_{\mathcal{H}_\epsilon(x)}\partial_\mathcal{S} \mathcal{W}^\epsilon(\mathcal{S},y-x)\frac{y-x}{|y-x|}\,dy,
\label{GradPD}
\end{eqnarray}
and
\begin{eqnarray}
\mathcal{S}=\frac{u^\epsilon(y)-u^\epsilon(x)}{|y-x|}\cdot e.
\label{sepsilon}
\end{eqnarray}
The displacement $u^\epsilon(t,x)$ is twice differentiable in time taking values in $L^2_0(D;\mathbb{R}^d)$. The space of such functions is denoted by $C^2([0,T];L^2_0(D;\mathbb{R}^d))$. The initial value problem for the peridynamic evolution  \eqref{stationary}
is seen to have a unique solution in this space, see Theorem \ref{existenceuniqueness} of Section \ref{EE}. The cohesive evolution $u^\epsilon(x,t)$ is uniformly bounded in the mean square norm over bounded time intervals $0<t<T$, i.e., 
\begin{eqnarray}
\max_{0<t<T}\left\{\Vert u^\epsilon(x,t)\Vert_{L^2(D;\mathbb{R}^d)}^2\right\}<C.
\label{bounds}
\end{eqnarray}
Here  $\Vert u^\epsilon(x,t)\Vert_{L^2(D;\mathbb{R}^d)}=\left(\int_D|u^\epsilon(x,t)|^2\,dx\right)^{1/2}$ and the upper bound $C$ is independent of $\epsilon$ and depends only on the initial conditions and body force applied up to time $T$, see Section \ref{GKP}.

The cohesive evolution has the following properties that are established in Section \ref{GKP}. The evolution has uniformly bounded kinetic and elastic potential energy
\begin{theorem}
\label{Gronwall}
{\rm \bf Bounds on kinetic and potential energy for cohesive dynamics}\\
There exists a positive constant $C$ depending only on $T$ and independent of   $\epsilon$ for which
\begin{eqnarray}
\sup_{0\leq t\leq T}\left\{PD^{\epsilon}(u^{\epsilon}(t))+\frac{\rho}{2}\Vert u_t^{\epsilon}(t)\Vert^2_{L^2(D;\mathbb{R}^d)}\right\}\leq C.
\label{boundenergy}
\end{eqnarray}
\end{theorem}

The evolution is uniformly continuous in time as measured by the mean square norm.
\begin{theorem}{\rm \bf Continuous cohesive evolution in mean square norm}\\
There is a positive constant $K$ independent of $t_2 < t_1$ in $[0,T]$ and  index $\epsilon$ for which
\begin{eqnarray}
\Vert u^{\epsilon}(t_1)-u^{\epsilon}(t_2)\Vert_{L^2(D;\mathbb{R}^d)}\leq K |t_1-t_2|.
\label{holderest}
\end{eqnarray}
\label{holdercont}
\end{theorem}

The evolution satisfies energy balance. The total energy of the cohesive evolution at time $t$ is given by
\begin{eqnarray}
\mathcal{EPD}^\epsilon(t,u^\epsilon(t))=\frac{\rho}{2}\Vert u^\epsilon_t(t)\Vert^2_{L^2(D;\mathbb{R}^d)}+PD^\epsilon(u^\epsilon(t))-\int_{D}b(t)\cdot u^\epsilon(t)\,dx
\label{energyt}
\end{eqnarray}
and the total energy of the system at time $t=0$ is
\begin{eqnarray}
\mathcal{EPD}^\epsilon(0,u^\epsilon(0))=\frac{\rho}{2}\Vert v_0\Vert^2_{L^2(D;\mathbb{R}^d)}+PD^\epsilon(u_0)-\int_{D}b(0)\cdot u_0\,dx.
\label{energy0}
\end{eqnarray}
The cohesive dynamics is seen to satisfy energy balance at every instant of the  evolution.
\begin{theorem}{\rm \bf Energy balance for cohesive dynamics}\\
\label{Ebalance}
\begin{eqnarray}
\mathcal{EPD}^\epsilon(t,u^\epsilon(t))=\mathcal{EPD}^\epsilon(0,u^\epsilon(0))-\int_0^t\int_{D} b_t(\tau)\cdot u^\epsilon(\tau)\,dx\,d\tau.\label{BalanceEnergy}
\end{eqnarray}
\end{theorem}

\setcounter{equation}{0} \setcounter{theorem}{0} \setcounter{lemma}{0}\setcounter{proposition}{0}\setcounter{remark}{0}\setcounter{remark}{0}\setcounter{definition}{0}\setcounter{hypothesis}{0}

\section{Dynamic instability and fracture nucleation}
\label{sec3}
In this section we present a fracture nucleation condition that arises from the unstable force law \eqref{forcestate}. This condition is manifested as a dynamic instability. In the following companion section we investigate the localization of dynamic instability as $\epsilon_k\rightarrow 0$ and define the notion of process zone for the cohesive evolution.
Fracture nucleation conditions can be viewed as instabilities and have been identified for peridynamic evolutions in \cite{SillingWecknerBobaru}. Fracture nucleation criteria formulated as instabilities for one dimensional peridynamic bars are developed in  \cite{WecknerAbeyaratne}. In this treatment we define a source for crack nucleation as  jump discontinuity in the displacement field that can become unstable and grow in time. Here we establish a direct link between the growth of jump discontinuities and the appearance of strain concentrations inside the deforming body.

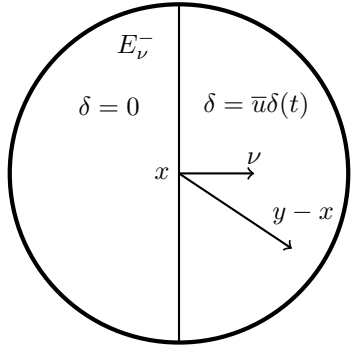
\begin{figure} 
\centering
\begin{tikzpicture}[xscale=1,yscale=1]
\draw [-,thick] (0,2.25) -- (0,-2.25);
\draw [->,thick] (0,0) -- (1,0);
\draw [->,thick] (0,0) -- (1.5,-1.0);
%\draw [-,thick] (-3,-0.25) to [out=-25,in=180] (-1.5,-1.5) to [out=0,in=180] (1.5,1.5)
%           to [out=0,in=165] (3,0.25);
%\draw (1.5,-0.2) -- (1.5, 0.2);
%\draw (-1.5,-0.2) -- (-1.5, 0.2);
\draw [ultra thick] (0,0) circle [radius=2.25];
\node [left] at (0,0) {$x$};
\node [above] at (1,0.0) {$\nu$};
\node [above] at (1.65,-0.80) {$y-x$};
%\node [below] at (-1.5,-0.2) {$-\frac{\overline{r}}{\sqrt{\epsilon|\xi|}}$};
%\node [right] at (3.5,0) {$\mathcal{S}$};
\node [right] at (0.2,0.9) {$\delta=\overline{u}\delta(t)$};
\node [left] at (-0.4,0.9) {$\delta=0$};
\node [right] at (-0.95,1.7) {$E^-_{\nu}$};
\end{tikzpicture} 
\caption{{\bf Jump discontinuity.}}
 \label{plane}
\end{figure}

We proceed with a formal perturbation analysis and consider a time independent body force density $b$ and a smooth equilibrium solution $u$ of \eqref{stationary}. Now perturb $u$ in the neighborhood of a point $x$ by adding a piecewise smooth discontinuity denoted by the vector field $\delta$. The perturbation takes the value zero on one side of a plane with normal vector $\nu$  passing through $x$ and on the other side of the plane takes the value $\delta=\overline{u}s(t)$. Here $s(t)$ is a scalar function of time and $\overline{u}$ is a constant vector. Consider the neighborhood $\mathcal{H}_\epsilon(x)$, then  $\delta(y)=0$ for $(y-x)\cdot\nu<0$ and $\delta(y)=\overline{u}s(t)$ for $(y-x)\cdot\nu\geq 0$, see Figure \ref{plane}. The half space on the side of the plane for which $(y-x)\cdot\nu<0$ is denoted by $E^-_{\nu}$.

Write  $u^p=u+\delta$ and assume
\begin{eqnarray}
\rho u^p_{tt}&=&-\nabla PD^\epsilon(u^p)+b.\label{stationarydiff1}
\end{eqnarray}
We regard $s(t)$ as a small perturbation and expand the integrand of $\nabla PD^\epsilon(u^p)$ in a Taylor series to recover the linearized evolution equation for the jump $s=s(t)$. The evolution equation is given by
\begin{eqnarray}
\rho s_{tt}\overline{u}=\mathcal{A}_{\nu}(x)\overline{u}s,
\label{pertevolution}
\end{eqnarray}
where the stability matrix $\mathcal{A}_{\nu}(x)$ is a $d\times d$ symmetric matrix with real eigenvalues and is defined by
\begin{eqnarray}
\mathcal{A}_{\nu}(x)&=&-\frac{2}{\epsilon V_d}\left\{\int_{\mathcal{H}_\epsilon(x)\cap E^-_{\nu}}\frac{1}{|y-x|}\partial^2_{\mathcal{S}}\mathcal{W}^\epsilon(\mathcal{S},y-x)\frac{y-x}{|y-x|}\otimes\frac{y-x}{|y-x|}dy\right\},
\label{instabilitymatrix}
\end{eqnarray}
and $$\mathcal{S}=\mathcal{S}(y,x)=\left(\frac{u(y)-u(x)}{|y-x|}\right)\cdot\frac{y-x}{|y-x|}.$$
Calculation shows that
\begin{eqnarray}
\partial^2_{\mathcal{S}}\mathcal{W}^\epsilon(\mathcal{S},y-x)=\frac{2}{\epsilon}J^\epsilon(|y-x|)\left(f'\left(|y|\mathcal{S}^2\right)+2f''\left(|y-x|\mathcal{S}^2\right)|y-x|\mathcal{S}^2\right),\label{expand}
\end{eqnarray}
where $f'(|y|\mathcal{S}^2)>0$ and $f''(|y|\mathcal{S}^2)<0$. On writing
\begin{eqnarray}
\mathcal{S}_c=\frac{\overline{r}}{\sqrt{|y-x|}}
\label{critstrain}
\end{eqnarray}
we have that
\begin{eqnarray}
\partial_\mathcal{S}^2\mathcal{W}^\epsilon(\mathcal{S},y)>0\hbox{    for   }|\mathcal{S}(y,x)|<\mathcal{S}_c,
\label{loss1}
\end{eqnarray} 
and 
\begin{eqnarray}
\partial_\mathcal{S}^2W^\epsilon(\mathcal{S},y)<0\hbox{    for   }|\mathcal{S}(y,x)|>\mathcal{S}_c. 
\label{loss2}
\end{eqnarray}
Here $\overline{r}$ is the inflection point for the function $r:\rightarrow f(r^2)$ and is the root of the equation
\begin{eqnarray}
f'(r)+2{r}f''({r})=0.
\label{rootroverline}
\end{eqnarray}
Note that the critical strain $\mathcal{S}_c$ for which the cohesive force between a pair of points $y$ and $x$ begins to soften is akin to the square root singularity seen at the crack tip in classical brittle fracture mechanics.

For eigenvectors $\overline{u}$ in the eigenspace associated with positive eigenvalues $\lambda$  of $\mathcal{A}_\nu(x)$ one has
\begin{eqnarray}
\rho\partial^2_{tt}s(t)=\lambda s(t)
\label{stabeq}
\end{eqnarray}
and the perturbation $s(t)$ can grow exponentially.  Observe from \eqref{loss2} that the quadratic form
\begin{eqnarray}
\mathcal{A}_{\nu}(x)\overline{w}\cdot\overline{w}=-\frac{2}{\epsilon V_d}\left\{\int_{\mathcal{H}_\epsilon(x)\cap E^-_{\nu}}\frac{1}{|y-x|}\partial^2_{\mathcal{S}}\mathcal{W}^\epsilon(\mathcal{S},y-x)(\frac{y-x}{|y-x|}\cdot\overline{w})^2dy\right\}
\label{quadformunstable}
\end{eqnarray}
will have at least one positive eigenvalue provided a sufficiently large proportion of bonds $y-x$ inside the horizon have strains satisfying 
\begin{eqnarray}
|\mathcal{S}(x,y)|>\mathcal{S}_c
\label{exceed}
\end{eqnarray} 
for which the cohesive force is in the unstable phase. For this case we see that the jump can grow exponentially.  The key feature here is that dynamic instability is explicitly linked to strain concentrations in this cohesive model as is seen from \eqref{loss2} together with \eqref{quadformunstable}. Collecting results we have the following proposition.

\begin{proposition}{\em \bf Fracture nucleation condition for  cohesive dynamics}\\
\label{Nuccriteria}
A condition for crack nucleation at a point $x$  is that there is at least one direction $\nu$ for which $\mathcal{A}_{\nu}(x)$ has at least one positive eigenvalue. This occurs if there is a square root strain concentration  $|\mathcal{S}(y,x)|>\mathcal{S}_c$ over a sufficiently large proportion of cohesive bonds inside the peridynamic horizon.
\end{proposition}
\noindent Proposition \ref{Nuccriteria} together with (\ref{loss2}) provide the explicit link between dynamic instability and  the critical  strain where the cohesive law begins to soften.

More generally  we may postulate a condition for  the  direction  along which the opposite faces of a nucleating fissure are oriented and the direction of the displacement jump across it. Recall that two symmetric matrices $A$ and $B$ satisfy $A\geq B$ in the sense of quadratic forms if $A\overline{w}\cdot\overline{w}\geq B\overline{w}\cdot\overline{w}$ for all $\overline{w}$ in $\mathbb{R}^d$. We say that a matrix $A$ is the maximum of a collection of symmetric matrices if $A\geq B$ for all matrices $B$ in the collection. 

We postulate that the faces of the nucleating fissure are perpendicular to the direction $\nu^*$ associated with the  the matrix $\mathcal{A}_{\nu^*}(x)$ for which
\begin{eqnarray}
\mathcal{A}_{\nu^*}(x)=\max\left\{\mathcal{A}_{\nu}(x);\, \hbox{over all directions } \nu\hbox{ such that } \mathcal{A}_{\nu}(x) \hbox{ has a positive eigenvalue}\right\},
\label{bestdirctionforgrowth}
\end{eqnarray}
and that the orientation of the jump in displacement across opposite sides of the fissure  lies in the eigenspace associated with the largest positive eigenvalue of $\mathcal{A}_{\nu^*}$,  {\em i.e., the fissure is oriented along the most unstable orientation and  the displacement jump across the nucleating fissure is along the most unstable direction.}

\setcounter{equation}{0} \setcounter{theorem}{0} \setcounter{lemma}{0}\setcounter{proposition}{0}\setcounter{remark}{0}\setcounter{remark}{0}\setcounter{definition}{0}\setcounter{hypothesis}{0}

\section{The process zone for cohesive dynamics and its localization in the small horizon limit }
\label{sec4}

In this section it is shown that the collection of centroids of peridynamic neighborhoods with strain exceeding a prescribed  threshold  concentrate on sets with zero volume in the limit of vanishing non-locality.   In what follows we probe the dynamics to obtain mathematically rigorous and explicit estimates on the size of the process zone in terms of the radius of the peridynamic horizon $0<\epsilon<1$.  

The continuity of the displacement inside the neighborhood $\mathcal{H}_\epsilon(x)$ is measured  quantitatively  by
\begin{eqnarray}
|u(y)-u(x)|\leq\underline{k}\,|y-x|^\alpha, \hbox{ for $y\in \mathcal{H}_\epsilon(x)$},
\label{moduluscont}
\end{eqnarray}
with $0<\underline{k}$ and  $0\leq\alpha\leq 1$. In what follows we focus on the reduction of continuity measured quantitatively by
\begin{eqnarray}
|u(y)-u(x)|>\underline{k}\,|y-x|^\alpha \hbox{ for $y\in \mathcal{H}_\epsilon(x)$}.
\label{modulusdiscont}
\end{eqnarray}
Observe when \eqref{modulusdiscont} holds and $\alpha=1/2$ and $\underline{k}=\overline{r}$ then $|\mathcal{S}(y,x)|>\mathcal{S}_c$  and there is softening in the cohesive force-strain behavior given by \eqref{forcestate} .

We now consider solutions $u^\epsilon$ of \eqref{stationary} and define a mathematical notion of process zone based the strain exceeding threshold values associated with \eqref{modulusdiscont}. 
The process zone is best described in terms of the basic unit of peridynamic interaction: the peridynamic neighborhoods $\mathcal{H}_\epsilon(x)$ of radius $\epsilon>0$ with centers $x\in D$.  We fix a choice of $\underline{k}$ and $\alpha$ belonging to the intervals $0< \underline{k}\leq \overline{r}$ and  $1/2\leq \alpha<1$.  The strain between $x$ and a point $y$ inside the neighborhood is denoted by $\mathcal{S}^\epsilon(y,x)$.
The collection of points $y$ inside $\mathcal{H}_\epsilon(x)$ for which the  strain $|\mathcal{S}^\epsilon(y,x)|$ exceeds the threshold function   $\underline{k}\,|y-x|^{\alpha-1}$ is denoted by
$\{y\hbox{ $ \in$ }\mathcal{H}_\epsilon(x):\, |\mathcal{S}^\epsilon(x,y)|>\underline{k}\,|y-x|^{\alpha-1}\}$. Note for $0<\underline{k}<\overline{r}$ and $1/2<\alpha<1$ that
\begin{eqnarray}
&&\left\{y\hbox{ $\in$ }\mathcal{H}_\epsilon(x):\, |\mathcal{S}^\epsilon(y,x)|>\mathcal{S}_c\right\}\subset\{y\hbox{ $\in$ }\mathcal{H}_\epsilon(x):\, |\mathcal{S}^\epsilon(y,x)|>\frac{\underline{k}}{|y-x|^{1-\alpha}}\}.\label{nonlipshitz}
\end{eqnarray}

%We introduce the indicator function $\chi^{+,\epsilon}_{\alpha,\underline{k}}(y,x)$ taking the value $1$ for $y$ belonging to $\{y\hbox{ $\in$ }%\mathcal{H}_\epsilon(x):\, |\mathcal{S}^\epsilon(y,x)|>\underline{k}\,|y-x|^{\alpha-1}\}$ and zero otherwise.
The fraction of points inside the neighborhood $\mathcal{H}_\epsilon(x)$ with strains exceed the threshold  is written
\begin{eqnarray}
P\left(\{y\hbox{ in }\mathcal{H}_\epsilon(x):\, |\mathcal{S}^\epsilon(y,x)|>\underline{k}\,|y-x|^{\alpha-1}\}\right),\label{weight}
\end{eqnarray}
where the weighted volume fraction for any subset $B$ of $\mathcal{H}_\epsilon(x)$ is defined as 
\begin{eqnarray}
P(B)=\frac{1}{\epsilon^d m}\int_{B}\,(|y-x|/\epsilon)J(\vert y-x\vert/\epsilon)\,dy,
\label{weightdefined}
\end{eqnarray}
with normalization constant 
\begin{eqnarray}
m=\int_{\mathcal{H}_1(0)}|\xi||J(|\xi|)\,d\xi
\label{normalize}
\end{eqnarray}
chosen so that $P(\mathcal{H}_\epsilon(x))=1$.
%The weighted volume fraction of strains inside $\mathcal{H}_\epsilon(x)$for which $ |\mathcal{S}^\epsilon(x,y)|>\mathcal{S}_c$ is bounded by
%\begin{eqnarray}
%P\left(\{y\hbox{ in }H_\epsilon(x):\, |\mathcal{S}^\epsilon(x,y)|>\mathcal{S}_c\}\right)\leq P\left(\{y\hbox{ in }H_\epsilon(x):\, |\mathcal{S}^\epsilon(y,x)|>\underline{k}\,|y-x|^{\alpha-1}\}\right),
%\label{unstablee}
%\end{eqnarray}
%for $0<\underline{k}\leq\overline{r}$ and $1/2\leq\alpha\leq 1$.

\begin{definition}{\bf Process Zone.}
\label{processZone}
Fix a volume fraction $0<\overline{\theta}\leq 1$, $0<\underline{k}\leq\overline{r}$, and $1/2\leq\alpha<1$ and at each time $t$ in the interval $0\leq t\leq T$, define the process zone $PZ^\epsilon(\underline{k},\alpha,\overline{\theta},t)$ to be the collection of centers of peridynamic neighborhoods for which the portion of points $y$ with strain $\mathcal{S}^\epsilon(t,y,x)$ exceeding the threshold $\underline{k}\,|y-x|^{\alpha-1}$ is greater than $\overline{\theta}$, i.e.,  $P\left(\{y\hbox{ in }\mathcal{H}_\epsilon(x):\, |\mathcal{S}^\epsilon(t,y,x)|>\underline{k}\,|y-x|^{\alpha-1}\}\right)>\overline{\theta}$. 
\end{definition}

The fracture set is defined to be the process zone for which strains exceed the  threshold $\mathcal{S}_c$ and the force vs. strain curve begins to soften.
\begin{definition}{\bf Fracture Set.}
\label{Fractureset}
The fracture set is defined to be the process zone associated with the values ${\theta}=1/2$, $\underline{k}=\overline{r}$, and $\alpha=1/2$ and at each time $t$ in the interval $0\leq t\leq T$, and is defined by $PZ^\epsilon(\overline{r},1/2,1/2,t)$ to be the collection of centers of peridynamic neighborhoods for which the portion of points $y$ with strain $\mathcal{S}^\epsilon(t,y,x)$ exceeding the threshold $\mathcal{S}_c$ is greater than $1/2$, i.e.,  $P\left(\{y\hbox{ in }\mathcal{H}_\epsilon(x):\, |\mathcal{S}^\epsilon(t,y,x)|>\mathcal{S}_c\}\right)>1/2$. 
\end{definition}
\noindent It is clear from the definition that the fracture set defined for this model contains the set of jump discontinuities for the displacement. This definition of fracture set given here is different that the usual one which collapses material damage onto a surface across which the displacement jumps.

It follows from Proposition \ref{Nuccriteria}  that the process zone contains peridynamic neighborhoods  associated with softening cohesive forces. Within this zone pre-existing jump discontinuities in the displacement field can grow.
\begin{remark}
Here we have described a range of process zones depending upon the choice of $\alpha$, $\underline{k}$ and $\overline{\theta}$. In what follows we show that for any choice of $\alpha$ in $1/2\leq \alpha <1$ and $\underline{k}$ in $0<\underline{k}\leq\overline{r}$ and $0<\overline{\theta}\leq 1$  the volume of the process zone is explicitly controlled by the radius of the peridynamic horizon $\epsilon$.
\end{remark}

We consider problem formulations in two and three dimensions and the volume or area of a set is given by the $d$ dimensional Lebesgue measure denoted by $\mathcal{L}^d$, for $d=2,3$.
We let 
\begin{eqnarray}
\label{upbdconst}
C(t)= \left( (2LEFM(u_0)+{\rho}\Vert v_0\Vert_{L^2(D;\mathbb{R}^d)}+1)^{1/2}+\sqrt{\rho^{-1}}\int_0^t\Vert b(\tau)\Vert_{L^2(D;\mathbb{R}^d)}\,d\tau\right)^2-1,
\end{eqnarray}
and note that $C(t)\leq C(T)$ for $t<T$.

We now give the following bound on the size of the process zone.
\begin{theorem} {\bf Dependence of the process zone on the radius of the peridynamic horizon}
\label{epsiloncontropprocesszone}
\begin{eqnarray}
\mathcal{L}^d\left(PZ^\epsilon(\underline{k},\alpha,\overline{\theta},t) \right)\leq \frac{\epsilon^{1-\beta}}{\overline{\theta}\underline{k}^2(f'(0)+o(\epsilon^\beta))}\times \frac{C(t)}{2m},
\label{controlbyepsilon}
\end{eqnarray}
where $0\leq\beta<1$ and $\beta=2\alpha-1$ and $0\leq t\leq T$.
\end{theorem}
\noindent Theorem \ref{epsiloncontropprocesszone} explicitly shows that the size of the process zone is controlled by the radius $\epsilon$ of the peridynamic horizon, uniformly in time. This theorem is proved in Section \ref{proofbondunstable}. 

\begin{remark}
\label{ModelParameter}
This analysis shows that {\em the horizon size $\epsilon$ for cohesive dynamics is a modeling parameter} that may be calibrated according to the size of the process zone obtained from experimental measurements.
\end{remark}

Next we show how the process zone localizes and concentrates on sets with zero volume in the small  horizon limit. To proceed choose  $\delta>0$  and consider  the sequence of solutions $u^{\epsilon_k}(t,x)$ to the cohesive dynamics for a family of radii  $\epsilon_k=\frac{1}{2^k}$, $k=1,\ldots$. The set of centers $x$ of neighborhoods $\mathcal{H}_{\epsilon_k}(x)$ that belong to  at least  one of the process zones $PZ^{\epsilon_k}(\underline{k},\alpha,\overline{\theta},t)$ for  some $\epsilon_k<\delta$ at time $t$ is denoted by $CPZ^\delta(\underline{k},\alpha,\overline{\theta},t)$. Let $CPZ^0(\underline{k},\alpha,\overline{\theta},t)=\cap_{0<\delta}CPZ^\delta(\underline{k},\alpha,\overline{\theta},t)$ be the collection  of centers of neighborhoods such that for every $\delta>0$ they belong to a process zone $PZ^{\epsilon_k}(\underline{k},\alpha,\overline{\theta},t)$ for  some $\epsilon_k<\delta$. The localization and concentration of the process zone is formulated in the following theorem.
\begin{theorem}{\rm\bf Localization of process zone in the small horizon limit.}\\
\label{bondunstable}
The collection of process zones  $CPZ^\delta(\underline{k},\alpha,\overline{\theta},t)$ is decreasing  with  $\delta\rightarrow 0$ and there is a positive constant $K$ independent of $t$ and $\delta$ for which
\begin{eqnarray}
&& \mathcal{L}^d\left(CPZ^\delta(\underline{k},\alpha,\overline{\theta},t)\right)\leq K{\delta^{1-\beta}}, \hbox{  for,  } 0\leq t\leq T,\,\,0<\beta=2\alpha-1\leq 1, \hbox{   with   }\nonumber\\
&& \mathcal{L}^d\left(CPZ^0(\underline{k},\alpha,\overline{\theta},t)\right)=\lim_{\delta\rightarrow 0}\mathcal{L}^d\left(CPZ^\delta(\underline{k},\alpha,\overline{\theta},t)\right)=0.
\label{limdelta}
\end{eqnarray}
For any choice of $0<\overline{\theta}\leq1$ the collection of centers of neighborhoods for which there exists a positive $\delta$  such that \begin{eqnarray}
P\left(\{y\hbox{ in }H_{\epsilon_k}(x):\, |\mathcal{S}^{\epsilon_k}(t,y,x)|\leq\underline{k}\,|y-x|^{\alpha-1}\}\right)\geq1-\overline{\theta},
\label{contolledstrain}
\end{eqnarray} 
for all $\epsilon_k<\delta$ is a set of full measure for every choice of $0<\underline{k}\leq\overline{r}$ and $1/2\leq\alpha<1$, i.e., $\mathcal{L}^d(D)=\mathcal{L}^d(D\setminus CPZ^0(\underline{k},\alpha,\overline{\theta},t))$.
\end{theorem}

\begin{remark}
\label{nearlipschitz}
The theorem shows that the process zone concentrates on a set of zero volume in the limit of vanishing peridynamic horizon. Note \eqref{contolledstrain} holds for any $0<\overline{\theta}\leq 1$. On  choosing $\overline{\theta}\approx 0$, and $\alpha\approx 1$ it is evident that the modulus of continuity for displacement field is close to Lipschitz outside of the process zone in the limit of vanishing nonolcality, $\epsilon_k\rightarrow 0$. The concentration of the process zone is inevitable for the cohesive model and is directly linked to the constraint on the energy budget  associated with the cohesive dynamics as described by Theorem \ref{Gronwall}.  This bound forces the localization of the process zone as shown in Section \ref{proofbondunstable}.
\end{remark}
\setcounter{equation}{0} \setcounter{theorem}{0} \setcounter{lemma}{0}\setcounter{proposition}{0}\setcounter{remark}{0}\setcounter{remark}{0}\setcounter{definition}{0}\setcounter{hypothesis}{0}

\section{The small horizon limit of cohesive dynamics}

\label{sec5}

In this section we identify the distinguished small horizon $\epsilon\rightarrow 0$ limit for cohesive dynamics.  It is shown here that the limit dynamics has bounded bulk linear elastic energy and Griffith surface energy characterized by the shear modului $\mu$, Lam\'e modulus $\lambda$,  and energy release rate $\mathcal{G}_c$ respectively. In order to make the connection between the limit dynamics and cohesive dynamics we will identify the relationship between the potentials $W^\epsilon(\mathcal{S},y-x)$ and the triple $\mu$, $\lambda$, $\mathcal{G}_c$. 

To reveal this connection consider a family  of cohesive evolutions $u^{\epsilon_k}$, each associated with a fixed potential $W^{\epsilon_k}$ and horizon length $\epsilon_k$, with $k=1,2,\ldots$ and $\epsilon_k\rightarrow 0$. Each $u^{\epsilon_k}(x,t)$ can be thought of as being  the result of a perfectly accurate numerical simulation of a cohesive evolution associated with the potential $W^{\epsilon_k}$. It is shown in this section  that the cohesive dynamics $u^{\epsilon_k}(x,t)$ converges to a limit evolution $u^0(x,t)$ in the limit, $\epsilon_k\rightarrow 0$. The limit evolution describes the dynamics of the cracked body when the scale of nonlocality is infinitesimally small with respect to the material specimen. Here the limiting free crack evolution is mediated through the triple $\mu$, $\lambda$, and $\mathcal{G}_c$ that are described by explicit formulas associated with the sequence of potentials $W^{\epsilon_k}$, see \eqref{calibrate1}, \eqref{calibrate2} and Theorem \ref{LEFMMThm} below. 

\bigskip

\noindent  {\em It is of fundamental importance to emphasize that we do do not impose a-priori relations between $W^{\epsilon_k}$ and the triple $\mu$, $\lambda$, and $\mathcal{G}_c$; we show instead that the cohesive dynamics $u^{\epsilon_k}(x,t)$  approaches the limit dynamics $u^0(x,t)$ characterized by $\mu$, $\lambda$, and $\mathcal{G}_c$ given by the formulas \eqref{calibrate1} and  \eqref{calibrate2} in the limit when $\epsilon_k\rightarrow 0$} .
 
 \bigskip
 
\noindent In what follows the  sequence of cohesive  dynamics described by $u^{\epsilon_k}$ is seen to converge to the limiting free crack evolution $u^0(x,t)$ in mean square, uniformly in time, see Theorem \ref{LimitflowThm}. 
The limit evolution is shown to have the following properties:
\begin{itemize}
\item It has  uniformly bounded  energy in the sense of linear elastic fracture mechanics for $0\leq t \leq T$.
\item It satisfies an energy inequality involving the kinetic energy of the motion together with the bulk elastic and surface energy associated with  linear elastic fracture mechanics for  $0\leq t\leq T$.
\end{itemize}

\noindent We  provide explicit conditions under which these properties are realized for the limit dynamics. 
\begin{hypothesis}
\label{remarkone}
We suppose that the magnitude  of the displacements $u^{\epsilon_k}$ for cohesive dynamics are bounded  for $0\leq t\leq T$ uniformly in $\epsilon_k$, i.e., $\sup_{\epsilon_k}\sup_{0\leq t\leq T}\Vert u^{\epsilon_k}(t)\Vert_{L^\infty(D;\mathbb{R}^d)}<\infty$.  
\end{hypothesis}

The convergence of  cohesive dynamics is given by the following theorem,
\begin{theorem} 
\label{LimitflowThm}
{\rm\bf Convergence of cohesive dynamics}\\
For each $\epsilon_k$ we prescribe identical LEFM initial data $u_0(x)$ and $v_0(x)$ and the solution to the cohesive dynamics initial value problem is denoted by $u^{\epsilon_k}$. Now consider a sequence of solutions $u^{\epsilon_k}$ associated with a vanishing peridynamic horizon $\epsilon_k\rightarrow 0$ and 
suppose Hypothesis  \ref{remarkone} holds true. Then, on passing to a subsequence if necessary, the cohesive evolutions $u^{\epsilon_k}$   converge in mean square uniformly in time to a limit evolution $u^0$ with the same LEFM initial data, i.e.,
\begin{eqnarray}
\lim_{\epsilon_k\rightarrow 0}\max_{0\leq t\leq T}\left\{\Vert u^{\epsilon_k}(t)-u^0(t)\Vert_{L^2(D;\mathbb{R}^d)}\right\}=0
\label{unifconvg}
\end{eqnarray}
and $u^0(x,0)=u_0(x)$ and $\partial_t u^0(x,0)=v_0(x)$.
\end{theorem}

To appropriately characterize the LEFM energy for the limit dynamics with freely propagating cracks one needs  a generalization of the strain tensor.  The appropriate notion of displacement and strain useful for problems involving discontinuities is provided by functions of bounded deformation BD introduced in \cite{Matthies}, \cite{Suquet}. The subspace of BD given by the special functions of bounded deformation $SBD$  introduced in \cite{Bellettini} is appropriate for describing discontinuities associated with linear elastic fracture. Functions in $u$  SBD belong to  $L^1(D;\mathbb{R}^d)$  and are approximately continuous, i.e., have Lebesgue limits for almost every $x\in D$  given by
\begin{eqnarray}
\lim_{r\searrow 0}\frac{1}{\pi r^d}\int_{B(x,r)}\,|u(y)-u(x)|\,dy=0, \hbox{ $d=2,3$}
\label{approx}
\end{eqnarray}
where $B(x,r)$ is the ball of radius $r$ centered at $x$.
The  jump set $J_{u}$  for elements of  $SBD$ is defined to be the set of points of discontinuity which have two different one sided Lebesgue limits.  One sided Lebesgue limits of  $u$ with respect to a direction $\nu_u(x)$ are  denoted by $u^-(x)$, $u^+(x)$ and are given by
\begin{eqnarray}
\lim_{r\searrow 0}\frac{1}{\pi r^d}\int_{B^-(x,r)}\,|u(y)-u^-(x)|\,dy=0,\,\,\, \lim_{r\searrow 0}\frac{1}{\pi r^d}\int_{B^+(x,r)}\,|u(y)-u^+(x)|\,dy=0, \hbox{$d=2,3$},
\label{approxjump}
\end{eqnarray}
where $B^-(x,r)$ and $B^+(x,r)$ are given by the intersection of $B(x,r)$ with the half spaces $(y-x)\cdot \nu_u(x)<0$ and $(y-x)\cdot \nu_u(x)>0$ respectively. SBD functions have jump sets $J_u$, described by a countable number of components $K_1,K_2,\ldots$, contained within smooth manifolds, with the exception of a set $K_0$ that has zero $d-1$ dimensional Hausdorff measure \cite{AmbrosioCosicaDalmaso}. Here the notion of arc length or (surface area)  is the $d-1$ dimensional Hausdorff  measure of $J_{u}$ and $\mathcal{H}^{d-1}(J_{u})=\sum_i\mathcal{H}^{d-1}(K_i)$.  
The  strain \cite{AmbrosioCosicaDalmaso} of a displacement $u$ belonging to SBD, written as $\mathcal{E}u$, is a generalization of the classic strain tensor and satisfies the property
\begin{eqnarray}
\lim_{r\searrow 0}\frac{1}{\pi r^d}\int_{B(x,r)}\,\frac{|(u(t,y)-u(t,x)-\mathcal{E}u(t,x)(y-x))\cdot(y-x)|}{|y-x|^2}\,dy=0, \hbox{ $d=2,3$}
\label{appgrad}
\end{eqnarray}
for almost every $x\in D$, with respect to $d$-dimensional Lebesgue measure $\mathcal{L}^d$. 
The symmetric part of the distributional derivative of $u$, $E u=1/2(\nabla u+\nabla u^T)$ for $SBD$ functions is a $d\times d$ matrix valued Radon measure with absolutely continuous part described by the density $\mathcal{E}u$ and singular part described by the jump set \cite{AmbrosioCosicaDalmaso}, \cite{Bellettini} and
\begin{eqnarray}
\langle E u,\Phi\rangle=\int_D\,\sum_{i,j=1}^d\mathcal{E}u_{ij}\Phi_{ij}\,dx+\int_{J_{u}}\,\sum_{i,j=1}^d(u^+_i - u^-_i)\nu_j\Phi_{ij}\,d\mathcal{H}^{d-1},
\label{distderiv}
\end{eqnarray}
for every continuous, symmetric matrix valued test function $\Phi$.
A description of $BD$ functions including their fine properties  and structure, together with the characterization of  $SBD$ functions on slices is developed in \cite{AmbrosioCosicaDalmaso} and \cite{Bellettini}.

The energy of linear elastic fracture mechanics extended to the class of $SBD$ functions is given by:
\begin{eqnarray}
LFEM(u,D)=\int_{D}\left(2\mu |\mathcal{E} u|^2+\lambda |{\rm div}\,u|^2\right)\,dx+\mathcal{G}_c\mathcal{H}^{d-1}(J_{u}), \hbox{ $d=2,3$,}
\label{LEFMSBVDefinition}
\end{eqnarray}
for $u$ belonging to $SBD$.
We now  describe the elastic energy for the limit dynamics. 
\begin{theorem}
{\rm\bf The limit dynamics has bounded LEFM energy}\\
The limit evolution $u^0$ belongs to SBD for every $t\in[0,T]$. Furthermore there exists a constant $C$ depending only on $T$  bounding the LEFM energy, i.e., 
\begin{eqnarray}
\int_{D}\,2\mu |\mathcal{E} u^0(t)|^2+\lambda |{\rm div}\,u^0(t)|^2\,dx+\mathcal{G}_c\mathcal{H}^{d-1}(J_{u^0(t)})\leq C, \hbox{ $d=2,3$,}
\label{LEFMbound}
\end{eqnarray}
for $0\leq t\leq T$. Here $\mu$, $\lambda$, and $\mathcal{G}_c$ are given by the explicit formulas
\begin{eqnarray}
\mu=\lambda=\frac{1}{4} f'(0)\int_{0}^1r^2J(r)dr, \hbox{ $d=2$}&\hbox{ and }& \mu=\lambda=\frac{1}{5} f'(0)\int_{0}^1r^3J(r)dr, \hbox{ $d=3$}\label{calibrate1}
\end{eqnarray}
and 
\begin{eqnarray}
\mathcal{G}_c=\frac{2\omega_{d-1}}{\omega_d}\, f_\infty \int_{0}^1r^dJ(r)dr, \hbox{  for $d=2.3$}
\label{calibrate2}
\end{eqnarray}
where $f_\infty$ is defined by \eqref{properties} and $\omega_{n}$ is the volume of the $n$ dimensional unit ball, $\omega_1=2,\omega_2=\pi,\omega_3=4\pi/3$.
The potential $f$ and influence function $J$ can always be chosen to satisfy \eqref{calibrate1} and \eqref{calibrate2} for any  $\mu=\lambda>0$ corresponding to the Poisson ratio $\nu=1/3$, for $d=2$ and $\nu=1/4$, for $d=3$, and $\mathcal{G}_c>0$.
\label{LEFMMThm}
\end{theorem} 

\begin{remark}
\label{boundedhausdorffmeasure}
The absolutely continuous part of the strain $\mathcal{E}u^0$ is defined for points away from the jump set $J_{u^0}$ and in this sense the process zone for the limit evolution can be viewed as being confined to the jump set $J_{u^0}$ .
Theorem \ref{LEFMMThm} shows that the jump set $J_{u^0}$ for the limit evolution $u^0(t,x)$ is confined to a set of finite $d-1$  dimensional Hausdorff measure.  
\end{remark}

We now present an energy inequality for the limit evolution. The sum of energy and work for the displacement $u^0$ at time $t$ is written
\begin{eqnarray}
\mathcal{GF}(u^0(t),D)=\frac{\rho}{2}\Vert u_t^0(t)\Vert^2_{L^2(D;\mathbb{R}^d)}+LEFM(u^0(t),D)-\int_{D}b(t)\cdot u^0(t)\,dx.
\label{sumtt}
\end{eqnarray}
The sum of energy and work for the initial data $u_0,v_0$ is written
\begin{eqnarray}
\mathcal{GF}(u_0,D)=\frac{\rho}{2}\Vert v_0\Vert^2_{L^2(D;\mathbb{R}^d)}+LEFM(u_0,D)-\int_{D}b(0)\cdot u_0\,dx.
\label{sumt0}
\end{eqnarray}
The energy inequality for the limit evolution $u^0$ is given by,
\begin{theorem} {\rm \bf Energy Inequality}\\
\label{energyinequality}
For almost every $t$ in $[0, T]$,
\begin{eqnarray}
\mathcal{GF}(u^0(t),D)\leq\mathcal{GF}(u_0,D)-\int_0^t\int_{D} b_t(\tau) \cdot u^0(\tau)\,dx\,d\tau.
\label{enegineq}
\end{eqnarray}
\end{theorem}

\begin{remark}
\label{remarkfinal}
The equality $\lambda=\mu$ appearing in  Theorem \ref {LEFMMThm}   is a consequence of the central force nature of the local cohesive interaction mediated by \eqref{forcestate}. More general non-central interactions are proposed in Section 15 of  \cite{Silling1} and in the state based peridynamic formulation  \cite{States}. The non-central formulations deliver a larger class of energy-volume-shape change relations for homogeneous deformations. 
Future work will address  state based formulations that deliver general anisotropic elastic response for the bulk energy associated with the limiting dynamics. 
\end{remark}

\setcounter{equation}{0} \setcounter{theorem}{0} \setcounter{lemma}{0}\setcounter{proposition}{0}\setcounter{remark}{0}\setcounter{remark}{0}\setcounter{definition}{0}\setcounter{hypothesis}{0}

\section{Free crack propagation in the small horizon limit}
\label{sec6}

We recall that the process zone concentrates on a set of zero volume (Lebesgue measure)  in the small horizon limit and identify conditions for which  the limit dynamics $u^0$ solves the  wave equation away from the evolving crack set. To begin we make a technical hypothesis on the regularity of the jump set of the limit dynamics $u^0(x,t)$.
\begin{hypothesis}
\label{remarkthree}
We  suppose that the crack set given by  $J_{u^0(t)}$ is a closed set for $0\leq t\leq T$.  
\end{hypothesis}
\noindent The next hypothesis applies to the concentration of the process zones as $\epsilon\rightarrow 0$ and their relation to the crack set for the limit dynamics.
\begin{hypothesis}
\label{remark555}
\noindent Theorem \ref{bondunstable} shows that the fracture sets defined as process zones with strains above $\mathcal{S}_c$, see Definition \ref{Fractureset},  concentrate on the set $CPZ^0(\overline{r},1/2,1/2,t)$. Here we assume that  $J_{u^0(t)}=CPZ^0(\overline{r},1/2,1/2,t)$  for   $0\leq t\leq T$.
\end{hypothesis}
\noindent The next hypothesis applies to neighborhoods $\mathcal{H}_{\epsilon_k}(x)$ for which the strain  is subcritical, i.e.,  $|\mathcal{S}^\epsilon|<\overline{r}/\sqrt{|y-x|}$, for $y$ in $\mathcal{H}_{\epsilon_k}(x)$. These neighborhoods will be referred to as neutrally stable.
\begin{hypothesis}
\label{remark556}
We suppose   that $\epsilon_k=\frac{1}{2^k}<\delta$ and $0\leq t\leq T$. Consider the collection of centers of peridynamic neighborhoods in $CPZ^\delta(\overline{r},1/2,1/2,t)$. We fatten out $CPZ^\delta(\overline{r},1/2,1/2,t)$ and consider $C\widetilde{PZ}^\delta(\overline{r},1/2,1/2,t)=\{x\in D:\, dist(x,CPZ^\delta(\overline{r},1/2,1/2,t)<\delta\}$. We suppose that all neighborhoods $H_{\epsilon_k}(x)$ that do not intersect the set $C\widetilde{PZ}^\delta(\overline{r},1/2,1/2,t)$ are neutrally stable.
\end{hypothesis}

\noindent  With these conditions satisfied the limit evolution $u^0$ is identified as a solution of the linear elastic wave equation.

\begin{theorem}
\label{waveequation}
Suppose Hypotheses \ref{remarkthree}, \ref{remark555} and \ref{remark556} hold true then the limit evolution $u^0(t,x)$ is a solution of the following wave equation (the first law of motion of Cauchy) in the sense of distributions on the domain $[0,T]\times D$ given by
\begin{eqnarray}
\rho u^0_{tt}= {\rm div}\sigma+b, \hbox{on $[0,T]\times D$},
\label{waveequationn}
\end{eqnarray}
where the stress tensor $\sigma$ is given by,
\begin{eqnarray}
\sigma =\lambda I_d Tr(\mathcal{E}\,u^0)+2\mu \mathcal{E}u^0,
\label{stress}
\end{eqnarray}
where $I_d$ is the identity on $\mathbb{R}^d$ and $Tr(\mathcal{E}\,u^0)$ is the trace of the strain. 
Here the second derivative $u_{tt}^0$ is the time derivative in the sense of distributions of $u^0_t$ and ${\rm div}\sigma$ is the divergence of the stress tensor $\sigma$ in the distributional sense.
\end{theorem}

\begin{remark}
\label{Disjointsets}
For completeness we recall that the strain $\mathcal{E} u^0(x,t)$ and jump set $J_{u^0(t)}$ are defined over disjoint sets in $[0,T]\times D$. 
\end{remark}

\begin{remark}
\label{displacementcrack}
The limit of the cohesive dynamics model is given by the displacement - crack set pair $u^0(t,x)$, $J_{u^0(t)}$. 
The wave equation provides the dynamic coupling between elastic waves and the evolving fracture path inside the media.
\end{remark}

\begin{remark}
\label{remarknearlyfinal}
Hypotheses \ref{remarkthree}, \ref{remark555}, and \ref{remark556}  are applied exclusively to establish Lemma  \ref{twolimitsB}  which identifies the absolutely continuous part of the limit strain
$\mathcal{E}_{ij}u^0e_ie_j$, $e=\xi/|\xi|$, with the weak $L^2 (D\times\mathcal{H}_1 (0))$ limit of the
strain $\mathcal{S}^\epsilon$restricted to pairs $(x,\xi)\in D\times\mathcal{H}_1 (0)$ for which $|\mathcal{S}^\epsilon|\leq\mathcal{S}_c$.
\end{remark}

\begin{remark}
\label{remarkfinal1}
We point out that the cohesive model addressed in this work does not have an irreversibility constraint and the constitutive law \eqref{forcestate} applies at all times in the peridynamic evolution.  Because of this the crack set at each time is given by $J_{u^0(t)}$. For rapid monotonic loading we anticipate that crack growth is increasing for this model, i.e., $J_{u^0(t')}\subset J_{u^0(t)}$ for $t'<t$. For cyclic loading this is clearly not the case and the effects of irreversibility (damage) must be incorporated into in the cohesive model. 
\end{remark}

\setcounter{equation}{0} \setcounter{theorem}{0} \setcounter{lemma}{0}\setcounter{proposition}{0}\setcounter{remark}{0}\setcounter{remark}{0}\setcounter{definition}{0}\setcounter{hypothesis}{0}

\section{Mathematical underpinnings and analysis}
\label{EightIntro}

In this section we provide the proofs of theorems stated in Sections \ref{sec2}, \ref{sec4}, \ref{sec5} and \ref{sec6}. The first subsection asserts the Lipschitz continuity of $\nabla PD^{\epsilon_k}(u)$ for $u$ in  $L^2_0(D;\mathbb{R}^d)$ and applies the   theory of ODE to deduce existence of the cohesive dynamics, see Section \ref{EE}. A Gronwall inequality is used to bound the cohesive potential energy and kinetic energy uniformly in time, see Section \ref{GKP}. Uniformly bounded sequences $\{u^{\epsilon_k}\}_{k=1}^\infty$ of cohesive dynamics  are shown to be compact in $C([0,T]; L^2_0(D;\mathbb{R}^d))$, see Section \ref{CC}. Any limit point $u^0$ for the sequence $u^{\epsilon_k}$ is shown to belong to SBD for every $0\leq t\leq T$, see Section \ref{CC}. The limit evolutions $u^0$ are shown to have uniformly bounded elastic energy in the sense of linear elastic fracture mechanics for $0\leq t\leq T$, see Section \ref{CC}. In Section \ref{EI} we pass to the limit in the energy balance equation for cohesive dynamics \eqref{BalanceEnergy} to recover an energy inequality for the limit flow. The wave equation satisfied by the limit flow is obtained on identifying the weak $L^2$ limit of the sequence $\{\nabla PD^{\epsilon_k}(u^{\epsilon_k})\}_{k=1}^\infty$  and passing to the limit in the weak formulation of \eqref{stationary}, see Section \ref{SC}. We conclude with the proof of Theorems \ref{epsiloncontropprocesszone} and  \ref{bondunstable}.

\subsection{Existence of a cohesive evolution}
\label{EE}

The peridynamic equation \eqref{eqofmotion} for cohesive dynamics is written as an equivalent first order system. We set $y^{\epsilon_k}=(y^{\epsilon_k}_1,y^{\epsilon_k}_2)^T$ where $y^{\epsilon_k}_1=u^{\epsilon_k}$ and $y_2^{\epsilon_k}=u_t^{\epsilon_k}$. Set $F^{\epsilon_k}(y^{\epsilon_k},t)=(F^{\epsilon_k}_1(y^{\epsilon_k},t),F^{\epsilon_k}_2(y^{\epsilon_k},t))^T$ where
\begin{eqnarray}
F^{\epsilon_k}_1(y^{\epsilon_k},t)&=&y_2^{\epsilon_k}\nonumber\\
F^{\epsilon_k}_2(y^{\epsilon_k},t)&=&\nabla PD^{\epsilon_k}(y_1^{\epsilon_k})+b(t).\nonumber
\end{eqnarray}
The initial value problem for $y^{\epsilon_k}$  given by the first order system is
\begin{eqnarray}
\frac{d}{dt} y^{\epsilon_k}=F^{\epsilon_k}(y^{\epsilon_k},t)\label{firstordersystem}
\end{eqnarray}
with initial conditions $y^{\epsilon_k}(0)=(u_0,v_0)^T$ satisfying LEFM initial conditions. In what follows we consider the more general class of initial data 
$(u_0,v_0)$ belonging to $L^2_0(D;\mathbb{R}^d)\times L^2_0(D;\mathbb{R}^d)$. 
\begin{theorem}
For $0\leq t\leq T$ there exists  unique solution in $C^1([0,T];L^2_0(D;\mathbb{R}^d)$ for the mesoscopic dynamics described by \eqref{firstordersystem} with initial data in $L^2_0(D;\mathbb{R}^d)\times L^2_0(D;\mathbb{R}^d)$ and body force $b(t,x)$ in $C^1([0,T];L^2_0(D;\mathbb{R}^d)$. 
\label{existenceuniqueness}
\end{theorem}
\noindent It now follows that for $LEFM$ initial data that one has a unique solution $u^{\epsilon_k}$ of \eqref{stationary}  in Section \ref{sec2}   belonging to $C^2([0,T];L^2_0(D;\mathbb{R}^d)$.

{\bf Proof of Theorem \ref{existenceuniqueness} }.
A straight forward calculation shows that for a generic positive constant $C$ independent of $\mathcal{S}$, $y-x$, and $\epsilon_k$, that
\begin{eqnarray}
\sup_{\mathcal{S}}|\partial_{\mathcal{S}}^2 W^{\epsilon_k}(\mathcal{S},y-x)|\leq \frac{C}{\epsilon_k|y-x|} \times J(|y-x|/\epsilon_k).
\label{secondderv}
\end{eqnarray}
From this it easily follows from H\"older and Minkowski inequalities that $\nabla PD^{\epsilon_k}$ is a Lipschitz continuous map from $L^2_0(D;\mathbb{R}^d)$ into $L^2_0(D;\mathbb{R}^d)$ and there is a positive constant $C$ independent of $0\leq t\leq T$, such that for any pair of vectors $y=(y_1,y_2)^T$, $z=(z_1,z_2)^T$ in $L^2_0(D;\mathbb{R}^d)\times L^2_0(D;\mathbb{R}^d)$ 
\begin{eqnarray}
\Vert F^{\epsilon_k}(y,t)-F^{\epsilon_k}(z,t)\Vert_{L^2(D;\mathbb{R}^d)^2}\leq \frac{C}{\epsilon_k}\Vert y-z\Vert_{L^2(D;\mathbb{R}^d)^2} \hbox{ for $0\leq t\leq T$}.
\label{lipschitz}
\end{eqnarray}
Here for any element $w=(w_1,w_2)$ of  $L^2_0(D;\mathbb{R}^d)\times L^2_0(D;\mathbb{R}^d)$, $\Vert w \Vert^2_{L^2(D;\mathbb{R}^d)^2}=\Vert w_1\Vert_{L^2(D;\mathbb{R}^d)}^2+\Vert w_2\Vert_{L^2(D;\mathbb{R}^d)}^2$.
Since  \eqref{lipschitz} holds the theory of ODE in Banach space \cite{Driver} shows that there exists a unique solution to the initial value problem \eqref{firstordersystem} with $y^{\epsilon_k}$ and $\partial_t y^{\epsilon_k}$ belonging to $C([0,T]; L^2_0(D;\mathbb{R}^d))$ and Theorem 
\ref{existenceuniqueness} is proved.
In this context we point out the recent work of \cite{EmmrichPhulst} where an existence theory for peridynamic evolutions for general pairwise force functions that are Lipschitz continuous with respect to the peridynamic deformation state is presented.

\subsection{Bounds on kinetic and potential energy for solutions of PD}
\label{GKP}
In this section we apply Gronwall's inequality to obtain bounds on the kinetic and elastic energy for peridynamic flows described by Theorem \ref{Gronwall}. The bounds are used to show that the solutions of the PD initial value problem are Lipschitz continuous in time

We now prove Theorem \ref{Gronwall}. Multiplying both sides of  \eqref{stationary} by $u_t^{\epsilon_k}(t)$ and integration together with a straight forward calculation gives
\begin{eqnarray}
&&\frac{1}{2}\frac{d}{dt}\left\{2PD^{\epsilon_k}(u^{\epsilon_k}(t))+{\rho}\Vert u_t^{\epsilon_k}(t)\Vert^2_{L^2(D;\mathbb{R}^d)}\right\}\nonumber\\
&&=\int_{D}(\nabla PD^{\epsilon_k}(u^{\epsilon_k}(t))+\rho u_{tt}^{\epsilon_k}(t))\cdot u_t^{\epsilon_k}(t)\,dx\nonumber\\
&&=\int_{D} u_t^{\epsilon_k}(t)\cdot b(t)\,dx\,\leq \, \Vert u_t^{\epsilon_k}\Vert_{L^2(D;\mathbb{R}^d)}\Vert b(t)\Vert_{L^2(D;\mathbb{R}^d)}.\label{esttime1}
\end{eqnarray}
Set
\begin{eqnarray}
&& W(t)=2PD^{\epsilon_k}(u^{\epsilon_k}(t))+{\rho}\Vert u_t^{\epsilon_k}(t)\Vert^2_{L^2(D;\mathbb{R}^d)}+1,
\label{wt}
\end{eqnarray}
and applying \eqref{esttime1} gives
\begin{eqnarray}
&&\frac{1}{2}W'(t) \leq \, \Vert u_t^{\epsilon_k}\Vert_{L^2(D;\mathbb{R}^d)}\Vert b(t)\Vert_{L^2(D;\mathbb{R}^d)}\leq\frac{1}{\sqrt{\rho}}\sqrt{W(t)}\Vert b(t)\Vert_{L^2(D;\mathbb{R}^d)}\label{esttime2}
\end{eqnarray}
and
\begin{eqnarray}
\frac{1}{2}\int_0^t\frac{W'(\tau)}{\sqrt{W(\tau)}}\,d\tau\leq\frac{1}{\sqrt{\rho}}\int_0^t\Vert b(\tau)\Vert_{L^2(D;\mathbb{R}^d)}\,d\tau.
\label{estime3}
\end{eqnarray}
Hence
\begin{eqnarray}
\sqrt{W(t)}-\sqrt{W(0)}\leq\frac{1}{\sqrt{\rho}}\int_0^t\Vert b(\tau)\Vert_{L^2(D;\mathbb{R}^d)}\,d\tau
\label{estime4}
\end{eqnarray}
and 
\begin{eqnarray}
&& 2PD^{\epsilon_k}(u^{\epsilon_k}(t))+{\rho}\Vert u_t^{\epsilon_k}(t)\Vert^2_{L^2(D;\mathbb{R}^d)} \leq \left(\frac{1}{\sqrt{\rho}}\int_0^t\Vert b(\tau)\Vert_{L^2(D;\mathbb{R}^d)}\,d\tau +\sqrt{W(0)}\right )^2-1.
\label{gineq}
\end{eqnarray}

For now we postpone the proof of \eqref{basicinequality} to Section \ref{CC} (see the discussion preceeding \eqref{upperboundperi}) and apply  \eqref{basicinequality} to get the upper bound 
\begin{eqnarray}
PD^{\epsilon_k}(u_0)\leq LEFM (u_0,D)\hbox{ for every $\epsilon_k$, \,\,$k=1,2,\ldots$},
\label{upperbound}
\end{eqnarray}
where $LEFM(u_0,D)$ is the elastic potential energy for linear elastic fracture mechanics given by \eqref{Gcrackenergy} or equivelently \eqref{LEFMSBVDefinition}. Theorem \ref{Gronwall} now follows from \eqref{gineq} and \eqref{upperbound}.

Theorem \ref{Gronwall} implies that PD solutions are Lipschitz continuous in time; this is stated explicitly in Theorem \ref{holdercont} of Section \ref{sec2}. To prove Theorem \ref{holdercont}  we write
\begin{eqnarray}
&&\Vert u^{\epsilon_k}(t_1)-u^{\epsilon_k}(t_2)\Vert_{L^2(D;\mathbb{R}^d)}=\left (\int_{D}|\int_{t_2}^{t_1} u^{\epsilon_k}_\tau(\tau)\,d\tau |^2\,dx\right )^{\frac{1}{2}}\nonumber\\
&&=\left (\int_{D}|t_1-t_2|^{2}\left|\frac{1}{|t_1-t_2|}\int_{t_2}^{t_1} u^{\epsilon_k}_\tau(\tau)\,d\tau \right|^2\,dx\right )^{\frac{1}{2}}\nonumber\\
&&\leq\left (\int_{D}|t_1-t_2|\int_{t_2}^{t_1} |u^{\epsilon_k}_\tau(\tau)|^2\,d\tau\,dx\right )^{\frac{1}{2}}\nonumber\\
&&\leq\left(|t_1-t_2|\int_{t_2}^{t_1}\Vert u_\tau^{\epsilon_k}(\tau)\Vert_{L^2(D;\mathbb{R}^d)}^2\,d\tau\right)^{1/2}\nonumber\\
&&\leq K|t_1-t_2|,
\label{lip}
\end{eqnarray}
where the third to last line follows from Jensen's inequality, the second to last line from Fubini's theorem and
the last inequality follows from the upper bound for $\Vert u_t^{\epsilon_k}(t)\Vert^2_{L^2(D;\mathbb{R}^d)}$ given by Theorem \ref{Gronwall}.

\subsection{Compactness and convergence}
\label{CC}
In this section we prove Theorems \ref{LimitflowThm} and \ref{LEFMMThm} . We start by establishing the inequality \eqref{basicinequality}  between the elastic energies $PD^{\epsilon_k}(u)$
and for $LEFM(u,D)$.  This is illustrated for any $u$ in $L^2_0(D;\mathbb{R}^d)\cap L^\infty(D;\mathbb{R}^d)$ and $LEFM(u,D)$ given by \eqref{LEFMSBVDefinition}. Here \eqref{LEFMSBVDefinition} reduces to \eqref{Gcrackenergy}  when $u$ is piecewise smooth and the crack $K$ consists of a finite number of smooth components. 
To obtain the upper bound we can directly apply the slicing technique of \cite{Gobbino3} to reduce to the one dimensional case to obtain an upper bound on one dimensional sections and then apply integral-geometric arguments to conclude.  Here the slicing theorem and integral-geometric measure appropriate for this approach in the context of SBD are given by Theorems  4.5 and  4.10 of  \cite{AmbrosioCosicaDalmaso}.
These arguments deliver the following inequality
\begin{eqnarray}
&&PD^{\epsilon_k}(u)\leq LEFM(u,D), \hbox{  for every $u$ in $L^2_0(D;\mathbb{R}^d)$, and $\epsilon_k>0$}. \label{upperboundperi}
\end{eqnarray}

To proceed with the proof of Theorem \ref{LimitflowThm} we require the compactness theorem.
\begin{theorem} {\rm \bf Compactness.}\\
Given a sequence of functions $u^{\epsilon_k}\in L^2_0(D;\mathbb{R}^d)$, $\epsilon_k=1/k,k=1,2,\ldots$ such that
\begin{eqnarray}
\sup_{\epsilon_k}\left(PD^{\epsilon_k}(u^{\epsilon_k})+\Vert u^{\epsilon_k}\Vert_{L^\infty(D;\mathbb{R}^d)}\right)<\infty,
\label{unifbound}
\end{eqnarray}
then there exists a subsequence $u^{{\epsilon_{k}}'}$ and limit point $u$  in $L_0^2(D;\mathbb{R}^d)$ for which
\begin{eqnarray}
u^{{\epsilon_k}'}\rightarrow u \hbox{ in $L^2(D;\mathbb{R}^d)$ as } {\epsilon_k}'\rightarrow 0.
\label{compactness}
\end{eqnarray}
\label{L2compact}
\end{theorem}

In what follows its convenient to change variables $y=x+\delta\xi$ for $|\xi|<1$ and $0<\delta<\alpha/2<1$, here the peridynamic neighborhood $\mathcal{H}_\delta(x)$ transforms to $\mathcal{H}_1(0)=\{\xi\in\mathbb{R}^d;\,|\xi|<1\}$. The unit vector $\xi/|\xi|$ is denoted by $e$.
To prove Theorem \ref{L2compact} we need the following upper bound given by the following theorem.
\begin{theorem}{\rm \bf Upper bound}\\
\label{coercivityth}
For any $0<\delta<\alpha/2$ there exists  positive constants $\tilde{K}_1$  and $\tilde{K}_2$ independent of $u\in L^2_0(D;\mathbb{R}^d)\cap L^\infty(D;\mathbb{R}^d)$ such that
\begin{eqnarray}
\int_{\mathcal{H}_1(0)}\int_D|u(x+\delta\xi)-u(x)|^2dx\,J(|\xi|)d\xi\leq \delta(\tilde{K}_1+\tilde{K}_2\Vert u\Vert_{L^\infty(D;\mathbb{R}^d)}^2) PD^\delta(u).
\label{coercivity}
\end{eqnarray}
\end{theorem}
We establish the upper bound in two steps.   
\begin{lemma}{\rm \bf Coercivity}\\
\label{coercivitya}
There exists a positive constant $C$ independent of $u\in L^2_0(D;\mathbb{R}^d)$
for which
\begin{eqnarray}
\int_{\mathcal{H}_1(0)}\int_D|u(x+\delta\xi)-u(x)|^2dx\,J(|\xi|)d\xi\leq C \int_{\mathcal{H}_1(0)}\int_D|(u(x+\delta\xi)-u(x))\cdot e|^2dx\,J(|\xi|)d\xi
\label{coercivea}
\end{eqnarray}
\end{lemma}

{\bf Proof of Lemma \ref{coercivitya}} We prove by contradiction. Suppose for every positive integer $N>0$ there is an element $u^N\in L_0^2(D;\mathbb{R}^d)$ for which
\begin{eqnarray}
&&N\int_{\mathcal{H}_1(0)}\int_D|(u^N(x+\delta\xi)-u^N(x))\cdot e |^2dx\,J(|\xi|)d\xi\nonumber\\
&&\leq\int_{\mathcal{H}_1(0)}\int_D|(u^N(x+\delta\xi)-u^N(x)) |^2dx\,J(|\xi|)d\xi.
\label{contra}
\end{eqnarray}
The Cauchy Schwartz inequality together with the triangle inequality deliver  a constant $\overline{K}>0$ for which
\begin{eqnarray}
\int_{\mathcal{H}_1(0)}\int_D|(u(x+\delta\xi)-u(x))|^2dx\,J(|\xi|)d\xi\leq\overline{K}\Vert u\Vert^2_{L^2(D;\mathbb{R}^d)}.
\label{upineq}
\end{eqnarray}
An application of the nonlocal Korn inequality, Lemma 5.5 of  \cite{DuGunzbergerlehoucqmengesha} gives the existence of a constant $\underline{K}>0$ independent of $u$ in $L_0^2(D;\mathbb{R}^d)$ for which
\begin{eqnarray}
\underline{K} \Vert u\Vert^2_{L^2(D;\mathbb{R}^d)}\leq\int_{\mathcal{H}_1(0)}\int_D|(u(x+\delta\xi)-u(x))\cdot e|^2dx\,J(|\xi|)d\xi.
\label{korn}
\end{eqnarray}
Applying the inequalities \eqref{contra}, \eqref{upineq}, and \eqref{korn} we discover that $\overline{K}/N\geq\underline{K}$ for all integers $N>0$ to conclude $\underline{K}=0$ which is a contradiction and Lemma \ref{coercivitya} is proved.
Theorem \ref{coercivityth} now follows from Lemma \ref{coercivitya} and the upper bound given by 
\begin{lemma} {\rm Upper bound}\\
\label{coercivityb}
\begin{eqnarray}
\int_{\mathcal{H}_1(0)}\int_D|(u(x+\delta\xi)-u(x))\cdot e|^2dx\,J(|\xi|)d\xi\leq \delta(\tilde{K}_1+\tilde{K}_2\Vert u\Vert_{L^\infty(D;\mathbb{R}^d)}^2) PD^\delta(u).
\label{coerciveb}
\end{eqnarray}
\end{lemma}
{\rm \bf Proof of Lemma \ref{coercivityb}}. Consider the concave potential function $f$ described in the introduction, recall $f(0)=0$ and given $M>0$ set $H_M=f(M)/M$. For $0<r<M$ one has $r<H_M^{-1}f(r)$ and set 
\begin{eqnarray}
A_{\delta\xi}=\{x\in D; |(u(x+\delta\xi)-u(x))\cdot e|^2>\delta|\xi|M\},
\label{Aset}
\end{eqnarray}
so
\begin{eqnarray}
&&\int_{D\setminus A_{\delta\xi}}\,|(u(x+\delta\xi)-u(x))\cdot e|^2\,dx=
\delta|\xi|\int_{D\setminus A_{\delta\xi}}\delta|\xi||\mathcal{S}|^2\,dx\nonumber\\
&&\leq\frac{\delta|\xi|}{H_M}\int_{D\setminus A_{\delta\xi}}\frac{1}{\delta|\xi|}f\left(\delta|\xi|\mathcal{S}|^2\right)\,dx.
\label{starter}
\end{eqnarray}
Now $f(r)>f(M)$ for $r>M$ gives
\begin{eqnarray}
\frac{1}{\delta|\xi|}f(M)\mathcal{L}^d(A_{\delta\xi})\leq
\int_{A_{\delta\xi}}\frac{1}{\delta|\xi|}f\left(\delta|\xi|\mathcal{S}|^2\right)\,dx 
\label{vol1}
\end{eqnarray}
and 
\begin{eqnarray}
\mathcal{L}^d(A_{\delta\xi}\label{ineqonA})\leq
\frac{\delta|\xi|}{f(M)}\int_{D}\frac{1}{\delta|\xi|}f\left(\delta|\xi|\mathcal{S}|^2\right)\,dx.
\label{volumebound}
\end{eqnarray}
Noting that
\begin{eqnarray}
\int_{A_{\delta\xi}}\,|(u(x+\delta\xi)-u(x))\cdot e|^2\,dx\leq 2\Vert u\Vert_{L^\infty(D;\mathbb{R}^d}^2\mathcal{L}^d(A_{\delta\xi})
\label{basiccc}
\end{eqnarray}
and collecting results one has
\begin{eqnarray}
\label{almostestimate}
\int_{D}\,|(u(x+\delta\xi)-u(x))\cdot e|^2\,dx\leq\delta|\xi|\left(\frac{1}{H_M}+\frac{2\Vert u\Vert_{L^\infty(D;\mathbb{R}^d)}^2}{f(M)}\right)
\int_{D}\frac{1}{\delta|\xi|}f\left(\delta|\xi|\mathcal{S}|^2\right)\,dx.
\end{eqnarray}
Lemma \ref{coercivityb} follows on multiplying both sides of \eqref{almostestimate} by $J(|\xi|)$ and integration over $\mathcal{H}_1(0)$.
Theorem \ref{coercivityth} follows from Lemmas \ref{coercivitya} and \ref{coercivityb}.

Arguing as in \cite{Gobbino3} we have the monotonicity given by
\begin{lemma}{\rm\bf Monotonicity}\\
For any integer $M$, $\eta>0$ and $u\in L^\infty(D;\mathbb{R}^d)$ one has
\begin{eqnarray}
PD^{M\eta}(u)\leq PD^\eta(u).
\label{monotone}
\end{eqnarray}
\label{monotonicity}
\end{lemma}
Now choose the subsequence $\epsilon_k=1/2^k,\,i=1,2,\ldots$ and from Theorem \ref{coercivityth} and Lemma \ref{monotonicity} we have for any $0<K<k$ with $\delta=2^{-K}$, $\epsilon_k=2^{-k}$,
\begin{eqnarray}
\int_{\mathcal{H}_1(0)}\int_D|u^{\epsilon_k}(x+\delta\xi)-u^{\epsilon_k}(x)|^2dx\,J(|\xi|)d\xi\leq \delta(\tilde{K}_1+\tilde{K}_2\Vert u^{\epsilon_k}\Vert_{L^\infty(D;\mathbb{R}^d)}^2) PD^{\epsilon_k}(u^{\epsilon_k}).\label{precompactness}
\end{eqnarray}
Applying the hypothesis \eqref{unifbound} to inequality \eqref{precompactness} gives a finite constant $B$ independent of $\epsilon_k$ and $\delta$ for which
\begin{eqnarray}
\int_{\mathcal{H}_1(0)}\int_D|u^{\epsilon_k}(x+\delta\xi)-u^{\epsilon_k}(x)|^2dx\,J(|\xi|)d\xi\leq \delta B,\label{precompactnesses}
\end{eqnarray}
for all $\epsilon_k<\delta$. One can then apply \eqref{precompactnesses} as in \cite{Gobbino3}, (or alternatively apply \eqref{precompactnesses} and arguments similar to the proof of the Kolomogorov-Riesz compactness theorem \cite{OlsenHolden}) to show that the sequence $\{u^{\epsilon_k}\}_{k=1}^\infty$ is a totally bounded subset of  $L^2_0(D;\mathbb{R}^d)$ and Theorem \ref{L2compact} is proved.

 Now it is shown that  the family of mesoscopic dynamics  $\{u^{\epsilon_k}\}_{k=1}^\infty$ is relatively compact in \\
$C([0,T];L_0^2(D;\mathbb{R}^d))$.
For each $t$  in $[0,T]$ we apply Theorem \ref{Gronwall} and Hypothesis \ref{remarkone} to obtain the bound
\begin{eqnarray}
PD^{\epsilon_k}(u^{\epsilon_k}(t))+\Vert u^{\epsilon_k}(t)\Vert_{L^\infty(D)}<C
\label{cpactbnd}
\end{eqnarray}
where $C<\infty$ and is independent of $\epsilon_k$, $k=1,2,\ldots$, and $0\leq t\leq T$.
With this bound we apply Theorem \ref{L2compact} to assert that for each $t$ the sequence $\{u^{\epsilon_k}(t)\}_{k=1}^\infty$
is relatively compact in $L^2(D;\mathbb{R}^d)$. 
From Theorem \ref{holdercont} the sequence $\{u^{\epsilon_k}\}_{k=1}^\infty$, is seen to be uniformly equi-continuous in $t$ with respect to the $L^2(D;\mathbb{R}^d)$ norm
and we immediately conclude from the Ascoli theorem that $\{u^{\epsilon_k}\}_{k=1}^\infty$ is relatively compact in $C([0,T];L^2(D;\mathbb{R}^d))$.
Therefore we can pass to a  subsequence also denoted by $\{u^{\epsilon_{k}}(t)\}_{k=1}^\infty$ to assert the existence of a limit evolution $u^0(t)$ in $C([0,T];L^2(D;\mathbb{R}^d))$  for which
\begin{eqnarray}
\lim_{k\rightarrow\infty}\left\{\sup_{t\in[0,T]}\Vert u^{\epsilon_{k}}(t)-u^0(t)\Vert_{\scriptscriptstyle{{L^2(D;\mathbb{R}^d)}}}\right\}=0
\label{unfconvergence}
\end{eqnarray}
and Theorem \ref{LimitflowThm} is proved.

We now prove Theorem \ref{LEFMMThm}. One has that limit points of sequences satisfying \eqref{unifbound}  enjoy higher regularity.

\begin{theorem}{\rm\bf Higher regularity}\\
\label{higherreg}
Every  limit point of a sequence $\{u^{\epsilon_k}\}_{k=1}^\infty$ in $L_0^2(D;\mathbb{R}^d)$ 
satisfying \eqref{unifbound}  belongs to $SBD$.
\end{theorem}

{\bf Proof.} To recover higher regularity one can directly apply the slicing technique of \cite{Gobbino3}  to reduce to the one dimensional case and construct sequences of functions converging in SBV to the limit point along one dimensional sections. One then applies  Theorem 4.7 of \cite{AmbrosioCosicaDalmaso} to conclude that the limit point belongs to $SBD$. 

It now follows from Theorem \ref{higherreg} that the limit evolution $u^0(t)$ belongs to $SBD$ for $0\leq t\leq T$.
Next we recall  the properties of $\Gamma$-convergence and apply them to finish the proof of Theorem \ref{LEFMMThm}. 
Consider a sequence of functions $\{F_j\}$ defined on a metric
space $\mathbb{M}$ with values in $\overline{\mathbb{R}}$ together with a function $F$ also defined on $\mathbb{M}$ with values in $\overline{\mathbb{R}}$.

\begin{definition}
\label{Gammaconvergence}
We say that $F$ is the $\Gamma$-limit of the sequence $\{F_j\}$ in $\mathbb{M}$   if the following two 
properties hold:
\begin{enumerate}
\item for every $x$ in $\mathbb{M}$ and every sequence $\{x_j\}$ converging to $x$, we have that
\begin{eqnarray}
F(x)\leq \liminf_{j\rightarrow\infty} F_j(x_j),\label{lowerbound}
\end{eqnarray}
\item for every $x$ in $\mathbb{M}$ there exists a recovery sequence $\{x_j\}$ converging to $x$, for which
\begin{eqnarray}
F(x)=\lim_{j\rightarrow\infty} F_j(x_j).\label{recovery}
\end{eqnarray}
\end{enumerate}
\end{definition}

For $u$ in $L^2_0(D;\mathbb{R}^d)$ define $PD^0:\,L^2_0(D;\mathbb{R}^d)\rightarrow [0,+\infty]$ by
\begin{equation}
PD^0(u)=\left\{ \begin{array}{ll}
LEFM(u,D)&\hbox{if $u$  belongs to $SBD$}\\
+\infty&\hbox{otherwise}
\end{array} \right.
\label{Gammalimit}
\end{equation}

A straight forward argument following  Theorem 4.3 $(ii)$ and $(iii)$ of \cite{Gobbino3} and invoking Theorems 4.5, 4.7, and 4.10 of \cite{AmbrosioCosicaDalmaso} as appropriate delivers

\begin{theorem}{\bf $\Gamma$- convergence and point wise convergence of peridynamic energies for cohesive dynamics.}
\begin{eqnarray}
&&PD^0 \hbox{ is the $\Gamma$-limit of $\{PD^{\epsilon_k}\}$ in $L^2_0(D;\mathbb{R}^d)$}, \hbox{ and }
\label{gammaconvpd}\\
&&\lim_{k\rightarrow\infty}PD^{\epsilon_k}(u)=PD^0(u), \hbox{  for every $u$ in $L^2_0(D;\mathbb{R}^d)$}.\label{pointwise}
\end{eqnarray}
\label{Gammaandpointwise}
\end{theorem}

Observe that since the sequence of peridynamic energies $\{PD^{\epsilon_k}\}$ $\Gamma$-converge to $PD^0$ in $L^2(D;\mathbb{R}^d)$ we can apply 
the lower bound property \eqref{lowerbound} of $\Gamma$-convergence to conclude that the limit has bounded elastic energy in the
sense of fracture mechanics, i.e.,
\begin{eqnarray}
LEFM(u^0(t),D)=PD^0(u^0(t))\leq\liminf_{k\rightarrow\infty}PD^{\epsilon_{k}}(u^{\epsilon_{k}}(t))<C.
\label{GSBV}
\end{eqnarray}
This concludes the proof of Theorem \ref{LEFMMThm}.

\subsection{Energy inequality for the limit flow}
\label{EI}
In this section we prove Theorem \ref{energyinequality}. We begin by showing that the limit evolution $u^0(t,x)$ has a weak derivative $u_t^0(t,x)$ belonging to $L^2([0,T]\times D;\mathbb{R}^d)$.  This is summarized in the following theorem.
\begin{theorem}
\label{weaktimederiviative}
On passage to subsequences if necessary the sequence $u_t^{\epsilon_k}$ weakly converges in $L^2([0,T]\times D;\mathbb{R}^d)$ to $u^0_t$ where
\begin{eqnarray}
-\int_0^T\int_D\partial_t\psi \cdot u^0\, dxdt=\int_0^T\int_D\psi \cdot u^0_t\, dxdt,
\label{weakl2time}
\end{eqnarray}
for all compactly supported smooth test functions $\psi$ on $[0,T]\times D$.
\end{theorem}

{\bf Proof.} The bound on the kinetic energy given in Theorem \ref{Gronwall} 
implies
\begin{eqnarray}
\sup_{\epsilon_k>0}\left(\sup_{0\leq t\leq T}\Vert u^{\epsilon_k}_t\Vert_{L^2(D;\mathbb{R}^d)}\right)< \infty.
\label{bddd}
\end{eqnarray}
Therefore the sequence $u^{\epsilon_k}_t$ is bounded in $L^2([0,T]\times D;\mathbb{R}^d)$ and passing to
a subsequence if necessary we conclude that there is a limit function
$\tilde{u}^0$ for which $u_t^{\epsilon_k}\rightharpoonup\tilde{u}^0$ weakly in $L^2([0,T]\times D;\mathbb{R}^d)$. Observe also that the uniform convergence \eqref{unfconvergence} implies that $u^{\epsilon_k}\rightarrow u^0$ in $L^2([0,T]\times D;\mathbb{R}^d)$. 
On writing the identity 
\begin{eqnarray}
-\int_0^T\int_D\partial_t\psi\cdot u^{\epsilon_k}\, dxdt=\int_0^T\int_D\psi \cdot u^{\epsilon_k}_t\, dxdt.
\label{weakidentity}
\end{eqnarray}
applying our observations and passing to the limit it is seen that
$\tilde{u}^0=u_t^0$ and the theorem follows.

To establish Theorem \ref{energyinequality} we require the following inequality.
\begin{lemma}
For almost every $t$ in $[0,T]$ we have
\label{weakinequality}
\begin{eqnarray}
\Vert u^0_t(t)\Vert_{L^2(D;\mathbb{R}^d)}\leq \liminf_{\epsilon_k\rightarrow 0}\Vert u^{\epsilon_k}_t(t)\Vert_{L^2(D;\mathbb{R}^d)}.
\label{limitweakineq}
\end{eqnarray}
\end{lemma}
{\bf Proof.}
We start with the identity
\begin{eqnarray}
\Vert u^{\epsilon_k}_t\Vert_{L^2(D;\mathbb{R}^d)}^2-2\int_D u^{\epsilon_k}_t\cdot u^0_t \,dx+\Vert u^0_t\Vert_{L^2(D;\mathbb{R}^d)}^2 
= \Vert u^{\epsilon_k}_t-u^0_t\Vert_{L^2(D;\mathbb{R}^d)}^2 \geq 0,
\label{locpositive}
\end{eqnarray}
and for every non-negative bounded measurable  function of time $\psi(t)$ defined on $[0,T]$ we have 
\begin{eqnarray}
\int_0^T\psi \Vert u^{\epsilon_k}_t-u^0_t\Vert_{L^2(D;\mathbb{R}^d)}^2\,dt\geq 0.
\label{positive}
\end{eqnarray}
Together with the weak convergence given in Theorem \ref{weaktimederiviative}  one easily sees that
\begin{eqnarray} 
\liminf_{\epsilon_k\rightarrow 0}\int_0^T\psi\Vert u^{\epsilon_k}_t\Vert_{L^2(D;\mathbb{R}^d)}^2\,dt-\int_0^T\psi\Vert u^0_t\Vert_{L^2(D;\mathbb{R}^d)}^2\,dt\geq 0.
\label{diff}
\end{eqnarray}

\noindent Applying  \eqref{bddd} and invoking the Lebesgue dominated convergence theorem we conclude
\begin{eqnarray}
\liminf_{\epsilon_k\rightarrow 0}\int_0^T\psi\Vert u^{\epsilon_k}_t\Vert_{L^2(D;\mathbb{R}^d)}^2\,dt=\int_0^T\psi\liminf_{\epsilon_k\rightarrow 0}\Vert u^{\epsilon_k}_t\Vert_{L^2(D;\mathbb{R}^d)}^2\,dt
\label{equalitybalance}
\end{eqnarray}
to recover the inequality given by
\begin{eqnarray}
\int_0^T\psi\left(\liminf_{\epsilon_k\rightarrow 0}\Vert u^{\epsilon_k}_t\Vert_{L^2(D;\mathbb{R}^d)}^2-\Vert u^0_t\Vert_{L^2(D;\mathbb{R}^d)}^2\right)\,dt\geq 0.
\label{diffinal}
\end{eqnarray}
The lemma follows noting that \eqref{diffinal} holds for every non-negative test function $\psi$.

Theorem \ref{energyinequality} now follows immediately on taking the $\epsilon_k\rightarrow 0$ limit in the peridynamic energy balance equation \eqref{BalanceEnergy} of Theorem \ref{Ebalance} and applying \eqref{pointwise},  \eqref{GSBV}, and \eqref{limitweakineq} of Lemma \ref{weakinequality}.

\subsection{Stationarity conditions for the limit flow}
\label{SC}

In this section we prove Theorem \ref{waveequation}. The first subsection establishes Theorem \ref{waveequation} using Theorem \ref{convgofelastic}. Theorem \ref{convgofelastic}  is proved in the second subsection.

\subsubsection{Proof of Theorem \ref{waveequation}}

To proceed we make the change of variables $y=x+\epsilon\xi$ where $\xi$ belongs to the unit disk $\mathcal{H}_1(0)$ centered at the origin
and the local strain $\mathcal{S}$ is of the form 
\begin{eqnarray}
\mathcal{S}=\left(\frac{u(x+\epsilon\xi)-u(x)}{\epsilon|\xi|}\right)\cdot e.
\label{strainrescaled}
\end{eqnarray}
It is convenient for calculation to express the strain through the directional difference operator $D_{e}^{\epsilon|\xi|}u$ defined by
\begin{eqnarray}
D_{e}^{\epsilon|\xi|}u(x)=\frac{u(x+\epsilon\xi)-u(x)}{\epsilon|\xi|}\hbox{ and }\mathcal{S}=D_{e}^{\epsilon|\xi|} u\cdot e,
\label{strainrescaledderiv}
\end{eqnarray}
with $e=\xi/|\xi|$. One also has
\begin{eqnarray}
D_{-e}^{\epsilon|\xi|}u(x)=\frac{u(x-\epsilon\xi)-u(x)}{\epsilon|\xi|},
\label{strainrescaledderivopposite}
\end{eqnarray}
and the integration by parts formula for functions $u$ in $L_0^2(D;\mathbb{R}^d)$, densities $\phi$ in $L_0^2(D;\mathbb{R})$ and $\psi$ continuous on $\mathcal{H}_1(0)$  given by
\begin{eqnarray}
\int_D\int_{\mathcal{H}_1(0)}(D_e^{\epsilon|\xi|}u\cdot e)\phi(x)\psi(\xi)\,d\xi\,dx=\int_D\int_{\mathcal{H}_1(0)}(u\cdot e) (D_{-e}^{\epsilon|\xi|}\phi)\psi(\xi)\,d\xi\,dx.
\label{intbypts}
\end{eqnarray}
Note further for $v$ in  $C^\infty_0(D;\mathbb{R}^d)$ and $\phi$ in $C^\infty_0(D;\mathbb{R})$ one has
\begin{eqnarray}
\lim_{\epsilon_k\rightarrow 0}D_e^{{\epsilon_k}|\xi|} v\cdot e=\mathcal{E} v \,e\cdot e &\hbox{ and }& \lim_{\epsilon_k\rightarrow 0}D_e^{{\epsilon_k}|\xi|} \phi= e\cdot \nabla\phi\label{grad}
\end{eqnarray}
where the convergence is uniform in $D$.

Taking the first variation of the action integral \eqref{Action} gives the Euler equation in weak form
\begin{eqnarray}
&&\rho\int_0^T\int_{D} u_t^{\epsilon_k}\cdot\delta_t\,dx \,dt\nonumber\\
&&-\frac{1}{\omega_d}\int_0^T\int_{D}\int_{\mathcal{H}_1(0)}|\xi|J(|\xi|)f'\left(\epsilon_k|\xi||D_e^{\epsilon_k|\xi|}u^{\epsilon_k}\cdot e|^2\right)2(D_e^{\epsilon_k|\xi|}u^{\epsilon_k}\cdot e) (D_{e}^{\epsilon_k|\xi|}\delta\cdot e)\,d\xi\,dx\,dt\nonumber\\
&&+\int_0^T\int_{D} b\cdot\delta\,dx\,dt=0,\label{stationtxweaker}
\end{eqnarray}
where the test function $\delta=\delta(x,t)=\psi(t)\phi(x)$ is smooth and has compact support in $[0,T]\times D$.
Next we make the change of function and write $F_s (\mathcal{S})=\frac{1}{s}f(s\mathcal{S}^2)$,  $F'_s(\mathcal{S})=2\mathcal{S}f'(s\mathcal{S}^2)$, and $s={\epsilon_k}|\xi|$ we transform \eqref{stationtxweaker} into
\begin{eqnarray}
&&\rho\int_0^T\int_{D} u_t^{\epsilon_k}\cdot\delta_t\,dx \,dt\nonumber\\
&&-\frac{1}{\omega_d}\int_0^T\int_{D}\int_{\mathcal{H}_1(0)}|\xi|J(|\xi|)F_{\epsilon_k|\xi|}'(D_{e}^{{\epsilon_k}|\xi|}u^{\epsilon_k}\cdot e)D_{e}^{\epsilon_k}\delta\cdot e\,d\xi\,dx\,dt\nonumber\\
&&+\int_0^T\int_{D} b\cdot\delta\,dx\,dt=0,\label{stationtxweakerlimitform}
\end{eqnarray}
where
\begin{eqnarray}
F_{\epsilon_k|\xi|}'(D_{e}^{\epsilon_k |\xi|}u^{\epsilon_k}\cdot e)=f'\left(\epsilon_k |\xi||D_{e}^{\epsilon_k |\xi|}u^{\epsilon_k}\cdot e|^2\right)2D_{e}^{\epsilon_k |\xi|}u^{\epsilon_k}\cdot e.
\label{deriv}
\end{eqnarray}

For future reference observe that $F_s(r)$ is convex-concave in $r$ with inflection point $\overline{r}_s=\overline{r}/\sqrt{s}$ where $\overline{r}$ is the inflection point of $f(r^2)=F_1(r)$. One also has the estimates
\begin{eqnarray}
&&F_s(r)\geq\frac{1}{s}F_1(\overline{r})\hbox{  for $r\geq\overline{r}_s$, and }\label{lowerestforF}\\
&&\sup_{0\leq r<\infty}|F'_s(r)|\leq\frac{2f'(\overline{r}^2){\overline{r}}}{\sqrt{s}},\label{boundderiv}
\label{estforFprime}
\end{eqnarray}
We send $\epsilon_k\rightarrow 0$ in \eqref{stationtxweakerlimitform} applying the weak convergence Theorem \ref{weaktimederiviative} to the first term to obtain
\begin{eqnarray}
&&\rho\int_0^T\int_{D} u_t^{0}\cdot\delta_t\,dx \,dt-\lim_{\epsilon_k\rightarrow 0}\frac{1}{\omega_d}\left(\int_0^T\int_{D}\int_{\mathcal{H}_1(0)}|\xi|J(|\xi|)F_{\epsilon_k|\xi|}'(D_{e}^{{\epsilon_k}|\xi|}u^{\epsilon_k}\cdot e)D_{e}^{\epsilon_k}\delta\cdot e\,d\xi\,dx\,dt\right)\nonumber\\
&&+\int_0^T\int_{D} b\cdot\delta\,dx\,dt=0.\label{stationtxweakerlimit}
\end{eqnarray}
Theorem \ref{waveequation} follows once we identify the limit of the second term in \eqref{stationtxweakerlimit} for smooth test functions $\phi(x)$ with support contained in $D$.
We state the following convergence theorem.
\begin{theorem}
\label{convgofelastic}
Given any infinitely differentiable test function $\phi$ with compact support in $D$ then
\begin{eqnarray}
\lim_{\epsilon_k\rightarrow 0}\frac{1}{\omega_d}\int_{D}\int_{\mathcal{H}_1(0)}|\xi|J(|\xi|)F_{\epsilon_k|\xi|}'(D_{e}^{{\epsilon_k}|\xi|}u^{\epsilon_k})D_{e}^{\epsilon_k}\phi\,d\xi\,dx=\int_{D}\mathbb{C}\mathcal{E} u^0:\mathcal{E}\phi\,dx,
\label{limitdelta}
\end{eqnarray}
where $\mathbb{C}\mathcal{E} u^0:\mathcal{E}\phi=\sum_{ijkl=1}^d\mathbb{C}_{ijkl}\mathcal{E}u^0_{ij}\mathcal{E}\phi_{kl}$, $\mathbb{C}\mathcal{E}u^0=\lambda I_d Tr(\mathcal{E}u^0)+2\mu \mathcal{E}u^0$, and  $\lambda$ and $\mu$ are given by \eqref{calibrate1}.
\end{theorem}
\noindent Theorem \ref{convgofelastic} is proved in Section \ref{prtheorem44}. 
The sequence of integrals on the left hand side of \eqref{limitdelta} are uniformly bounded in time, i.e.,
\begin{eqnarray}
\sup_{\epsilon_k>0}\left\{\sup_{0\leq t\leq T}\left\vert\int_{D}\int_{\mathcal{H}_1(0)}|\xi|J(|\xi|)F_{\epsilon_k|\xi|}'(D_{e}^{{\epsilon_k}|\xi|}u^{\epsilon_k}\cdot e)D_{e}^{\epsilon_k}\phi\cdot e\,d\xi\,dx\right\vert\right\}<\infty,
\label{uniboundt}
\end{eqnarray}
this is demonstrated in \eqref{Fprimesecond} of Lemma \ref{estimates}  in Section \ref{prtheorem44}. Applying the Lebesgue bounded convergence theorem together with Theorem \ref{convgofelastic}
with $\delta(t,x)=\psi(t)\phi(x)$ delivers the desired result
\begin{eqnarray}
&&\lim_{\epsilon_k\rightarrow 0}\frac{1}{\omega_d}\left(\int_0^T\int_{D}\int_{\mathcal{H}_1(0)}|\xi|J(|\xi|)F_{\epsilon_k|\xi|}'(D_{e}^{{\epsilon_k}|\xi|}u^{\epsilon_k}\cdot e)\psi D_{e}^{\epsilon_k}\phi\cdot e\,d\xi\,dx\,dt\right)\nonumber\\
&&=\int_0^T\int_{D}\mathbb{C}\mathcal{E} u^0:\mathcal{E}\phi\,dx\,dt,
\label{limitidentity}
\end{eqnarray}
and we recover the identity
\begin{eqnarray}
&&\rho\int_0^T\int_{D} u_t^{0}(t,x)\cdot\psi_t(t)\phi(x)\,dx \,dt-\int_0^T\int_{D}\psi(t)\mathbb{C}\mathcal{E} u^0(t,x):\mathcal{E}\phi(x)\,dx\,dt
\nonumber\\
&&+\int_0^T\int_{D} b(t,x)\cdot\psi(t)\phi(x)\,dx\,dt=0
\label{finalweakidentity}
\end{eqnarray}
from which Theorem \ref{waveequation} follows.

\subsubsection{ Proof of Theorem \ref{convgofelastic}}
\label{prtheorem44}

We decompose the difference $D_e^{\epsilon_k |\xi|}u^{\epsilon_k}\cdot e$ as 
\begin{eqnarray}
D_e^{\epsilon_k |\xi|}u^{\epsilon_k}\cdot e=(D_e^{\epsilon_k |\xi|}u^{\epsilon_k}\cdot e)^- +(D_e^{\epsilon_k |\xi|}u^{\epsilon_k}\cdot e)^+
\label{decompose}
\end{eqnarray}
where
\begin{equation}
(D_e^{\epsilon_k |\xi|}u^{\epsilon_k}(x)\cdot e)^-=\left\{\begin{array}{ll}(D_e^{\epsilon_k |\xi|}u^{\epsilon_k}(x)\cdot e),&\hbox{if  $|D_e^{\epsilon_k |\xi|}u^{\epsilon_k}(x)\cdot e|<\frac{\overline{r}}{\sqrt{\epsilon_k|\xi|}}$}\\
0,& \hbox{otherwise}
\end{array}\right.
\label{decomposedetails}
\end{equation}
where $\overline{r}$ is the inflection point for the function $F_1(r)=f(r^2)$. Here $(D_e^{\epsilon_k |\xi|}u^{\epsilon_k}\cdot e)^+$ is defined so that \eqref{decompose} holds.
We prove Theorem \ref{convgofelastic} by using the following two identities
described in the Lemmas below.
\begin{lemma}
\label{twolimitsA}
For any $\phi$ in $C^\infty_0(D;\mathbb{R}^d)$ 
\begin{eqnarray}
&&\lim_{\epsilon_k\rightarrow 0} \frac{1}{\omega_d}\int_{D}\int_{\mathcal{H}_1(0)}|\xi|J(|\xi|)F_{\epsilon_k|\xi|}'(D_{e}^{{\epsilon_k}|\xi|}u^{\epsilon_k}\cdot e)D_{e}^{\epsilon_k|\xi|}\phi\cdot e\,d\xi\,dx\nonumber\\
&&-2\lim_{\epsilon_k\rightarrow 0}\frac{1}{\omega_d}\int_{D}\int_{\mathcal{H}_1(0)}|\xi|J(|\xi|)f'(0)(D_e^{\epsilon_k |\xi|}u^{\epsilon_k}\cdot e)^-D_e^{\epsilon_k |\xi|}\phi\cdot e\,d\xi\,dx=0.
\label{functiotolimit}
\end{eqnarray}
\end{lemma}
\begin{lemma}
\label{twolimitsB}
Assume that Hypotheses \ref{remarkthree}, \ref{remark555} and \ref{remark556} hold true and
define the weighted Lebesgue measure $\nu$ by $d\nu=|\xi|J(|\xi|)d\xi\,dx$ for any Lebesgue measurable set $S\subset D\times\mathcal{H}_1(0)$.
Passing to subsequences if necessary  $\{(D_e^{\epsilon_k |\xi|}u^{\epsilon_k}\cdot e)^-\}_{k=1}^\infty$ converges weakly in $L^2(D\times\mathcal{H}_1(0);\nu)$ to $\mathcal{E} u^0 e\cdot e$, i.e., 
\begin{eqnarray}
&&\lim_{\epsilon_k\rightarrow 0}\frac{1}{\omega_d}\int_{D}\int_{\mathcal{H}_1(0)}(D_e^{\epsilon_k |\xi|}u^{\epsilon_k}\cdot e)^- \psi\,d\nu\nonumber\\
&&=\frac{1}{\omega_d}\int_{D}\int_{\mathcal{H}_1(0)}(\mathcal{E} u^0 e \cdot e)\psi\,d\nu,
\label{functiotofunction}
\end{eqnarray}
for any test function $\psi(x,\xi)$ in  $L^2(D\times\mathcal{H}_1(0);\nu)$. 
\end{lemma}

We now apply the Lemmas. Observing that $D_e^{\epsilon_k |\xi|}\phi\cdot e $ converges strongly in $L^2(D\times\mathcal{H}_1(0):\nu)$ to $\mathcal{E}\phi\, e\cdot e$ for test functions $\phi$ in $C^\infty_0(D;\mathbb{R}^d)$ and from the weak $L^2(D\times\mathcal{H}_1(0):\nu)$ convergence of $(D_e^{\epsilon_k |\xi|}u^{\epsilon_k}\cdot e)^-$ we deduce that
\begin{eqnarray}
&&\lim_{\epsilon_k\rightarrow 0}\frac{1}{\omega_d}\int_{D}\int_{\mathcal{H}_1(0)}|\xi|J(|\xi|)f'(0)(D_e^{\epsilon_k |\xi|}u^{\epsilon_k}\cdot e)^-(D_e^{\epsilon_k |\xi|}\phi\cdot e)\,d\nu\nonumber\\
&&=\frac{1}{\omega_d}\int_{D}\int_{\mathcal{H}_1(0)}|\xi|J(|\xi|)f'(0)\,(\mathcal{E} u^0 e\cdot e)(\mathcal{E}\phi e\cdot e)\,d\nu\nonumber\\
&&=\frac{f'(0)}{\omega_d}\sum_{ijkl=1}^d\int_{\mathcal{H}_1(0)}|\xi|J(|\xi|)\,e_i e_j e_k e_l,d\xi\int_{D}\mathcal{E}u^0_{ij}\mathcal{E}\phi_{kl}\,dx.\label{limitproduct}
\end{eqnarray}
Now we show that
\begin{eqnarray}
\frac{f'(0)}{\omega_d}\int_{\mathcal{H}_1(0)}|\xi|J(|\xi|)\,e_i e_j e_k e_l\,d\xi=\mathbb{C}_{ijkl}=2\mu\left(\frac{\delta_{ik}\delta_{jl}+\delta_{il}\delta_{jk}}{2}\right)+\lambda \delta_{ij}\delta_{kl}\label{shear}
\end{eqnarray}
where $\mu$ and $\lambda$ are given by \eqref{calibrate1}. To see this we write 
\begin{eqnarray}
\Gamma_{ijkl}(e)=e_ie_je_ke_l,
\label{tensorId}
\end{eqnarray}
to observe that $\Gamma(e)$ is a totally symmetric tensor valued function defined for  $e\in S^{d-1}$ with the property
\begin{eqnarray}
\Gamma_{ijkl}(Qe)=Q_{im}e_mQ_{jn}e_nQ_{ko}e_oQ_{lp}e_p=Q_{im}Q_{jn}Q_{ko}Q_{lp}\Gamma_{mnop}(e)
\label{rot}
\end{eqnarray}
for every rotation $Q$ in $SO^d$. Here repeated indices indicate summation.
We write
\begin{eqnarray}
\int_{\mathcal{H}_1(0)}|\xi|J(|\xi|)\,e_i e_j e_k e_l\,d\xi=\int_0^1|\xi|d|\xi|\int_{S^{d-1}}\Gamma_{ijkl}(e)\,de
\label{groupavgpre}
\end{eqnarray}
to see that for every $Q$ in $SO^d$
\begin{eqnarray}
Q_{im}Q_{jn}Q_{ko}Q_{lp}\int_{S^{d-1}}\Gamma_{ijkl}(e)\,de=\int_{S^{d-1}}\Gamma_{mnop}(Qe)\,de=\int_{S^{d-1}}\Gamma_{mnop}(e)\,de.
\label{groupavg}
\end{eqnarray}
Therefore we conclude that $\int_{S^{d-1}}\Gamma_{ijkl}(e)\,de$ is an isotropic symmetric $4^{th}$ order  tensor and of the form
\begin{eqnarray}
\int_{S^{d-1}}\Gamma_{ijkl}(e)\,de=a\left(\frac{\delta_{ik}\delta_{jl}+\delta_{il}\delta_{jk}}{2}\right)+b \delta_{ij}\delta_{kl}.
\label{groupavgaft}
\end{eqnarray}
Here we evaluate $a$  by contracting both sides of \eqref{groupavgaft} with a trace free matrix and $b$ by contracting both sides with the $d\times d$ identity and calculation delivers \eqref{shear}.
Theorem \ref{convgofelastic} now follows immediately from \eqref{limitproduct} and \eqref{functiotolimit}. 

To establish Lemmas \ref{twolimitsA} and \ref{twolimitsB}  we develop the following estimates for the sequences\\ $(D_e^{\epsilon_k |\xi|}u^{\epsilon_k}\cdot e)^-$ and $(D_e^{\epsilon_k |\xi|}u^{\epsilon_k}\cdot e)^+$.  We define the set $K^{+,\epsilon_k}$ by
\begin{eqnarray}
K^{+,\epsilon_k}=\{(x,\xi)\in D\times\mathcal{H}_1(0)\,: (D_e^{\epsilon_k |\xi|}u^{\epsilon_k}\cdot e)^+\not=0\}.
\label{suppset}
\end{eqnarray}

We have the following string of estimates.
\begin{lemma}
We introduce the generic positive constant $0<C<\infty$ independent of $0<\epsilon_k<1$ and $0\leq t\leq T$ and state the following inequalities that hold for all $0<\epsilon_k<1$ and $0\leq t\leq T$ and for $C^\infty(D)$ test functions $\phi$  with compact support on $D$.
\label{estimates}
\begin{eqnarray}
&&\int_{K^{+,\epsilon_k}} |\xi|J(|\xi|)\,d\xi\,dx<C\epsilon_k,\label{suppsetupper}\\
&&\left\vert\int_{D\times\mathcal{H}_1(0)}|\xi| J(|\xi|)F_{\epsilon_k|\xi|}'((D_{e}^{{\epsilon_k}|\xi|,+}u^{\epsilon_k}\cdot e)^+)(D_{e}^{\epsilon_k}\phi\cdot e)\,d\xi\,dx\right\vert<C\sqrt{\epsilon_k}\Vert\mathcal{E}\phi\Vert_{\scriptscriptstyle{{L^\infty(D;\mathbb{R}^{d\times d})}}},\label{Fprimefirst}\\
&&\int_{D\times\mathcal{H}_1(0)}|\xi| J(|\xi|)|(D_{e}^{{\epsilon_k}|\xi|}u^{\epsilon_k}\cdot e)^-|^2\,d\xi\,dx<C,\label{L2bound}\\
&&\int_{D\times\mathcal{H}_1(0)}|\xi| J(|\xi|)|D_{e}^{{\epsilon_k}|\xi|}u^{\epsilon_k}\cdot e|\,d\xi\,dx<C,\hbox{  and}\label{L1bound}\\
&&\left\vert\int_{D\times\mathcal{H}_1(0)}|\xi| J(|\xi|)F_{\epsilon_k|\xi|}'(D_{e}^{{\epsilon_k}|\xi|}u^{\epsilon_k}\cdot e) (D_{e}^{\epsilon_k}\phi\cdot e)d\,\xi\,dx\right\vert<C\Vert\mathcal{E} \phi\Vert_{\scriptscriptstyle{L^\infty(D;\mathbb{R}^{d\times d})}}.\label{Fprimesecond}
\end{eqnarray}
\end{lemma}
{\bf Proof.}
For $(x,\xi)\in K^{+,\epsilon_k}$ we apply \eqref{lowerestforF} to get
\begin{eqnarray}
J(|\xi|)\frac{1}{\epsilon_k}F_{1}(\overline{r})=|\xi|J(|\xi|)\frac{1}{\epsilon_k|\xi|}F_{1}(\overline{r})
\leq|\xi|J(|\xi|)F_{\epsilon_k|\xi|}(D_{e}^{{\epsilon_k}|\xi|}u^{\epsilon_k}\cdot e)\label{first}
\end{eqnarray}
and in addition since $|\xi|\leq 1$ we have
\begin{eqnarray}
&&\frac{1}{\epsilon_k}F_{1}(\overline{r}) \int_{K^{+,\epsilon_k}} |\xi|J(|\xi|)\,d\xi\,dx\leq\frac{1}{\epsilon_k}F_{1}(\overline{r}) \int_{K^{+,\epsilon_k}} J(|\xi|)\,d\xi\,dx\nonumber\\
&&\leq\int_{K^{+,\epsilon_k}} |\xi|J(|\xi|)F_{\epsilon_k|\xi|}(D_{e}^{{\epsilon_k}|\xi|}u^{\epsilon_k}\cdot e)\,d\xi\,dx\leq\sup_{t\in [0,T]}\sup_{\epsilon_k}PD^{\epsilon_k}(u^{\epsilon_k})
\label{firstb}
\end{eqnarray}
where Theorem \ref{Gronwall} implies that the right most element of  the sequence of inequalities is bounded and \eqref{suppsetupper} follows noting that the inequality \eqref{firstb} is equivalent to \eqref{suppsetupper}. More generally since $|\xi|\leq 1$ we may argue as above to conclude that 
\begin{eqnarray}
\int_{K^{+,\epsilon_k}} |\xi|^pJ(|\xi|)\,d\xi\,dx<C\epsilon_k.
\label{power}
\end{eqnarray}
for $0\leq p$.
We apply \eqref{estforFprime} and \eqref{power} to find
\begin{eqnarray}
&&\left\vert\int_{D\times\mathcal{H}_1(0)}|\xi| J(|\xi|)F_{\epsilon_k|\xi|}'((D_{e}^{{\epsilon_k}|\xi|}u^{\epsilon_k}\cdot e)^+)D_{e}^{\epsilon_k}\phi\,d\xi\,dx\right\vert\nonumber\\
&&\leq C\frac{2f'(\overline{r}^2)\overline{r}}{\sqrt{\epsilon_k}}\int_{K^{+,\epsilon_k}}\sqrt{|\xi|}J(|\xi|)\,d\xi\,dx \Vert\mathcal{E}\phi\Vert_{\scriptscriptstyle{L^\infty(D;\mathbb{R}^{d\times d})}}
\leq\sqrt{\epsilon_k}C\Vert\mathcal{E}\phi\Vert_{\scriptscriptstyle{L^\infty(D;\mathbb{R}^{d\times d})}},
\label{second}
\end{eqnarray}
and \eqref{Fprimefirst} follows.

A basic calculation shows there exists a positive constant independent of $r$ and $s$ for which
\begin{eqnarray}
r^2\leq C F_s(r), \hbox{  for $r<\frac{\overline{r}}{\sqrt{s}}$},
\label{rsquaredd}
\end{eqnarray}
so
\begin{eqnarray}
|D_{e}^{\epsilon_k |\xi|}u^{\epsilon_k}\cdot e|^2\leq C F_{\epsilon_k|\xi|}(D_{e}^{{\epsilon_k}|\xi|}u^{\epsilon_k}\cdot e), \hbox{  for $|D_{e}^{\epsilon_k|\xi|}u^{\epsilon_k}\cdot e|<\frac{\overline{r}}{\sqrt{\epsilon_k |\xi|}}$},
\label{rsquared}
\end{eqnarray}
and
\begin{eqnarray}
&&\int_{D\times\mathcal{H}_1(0)}|\xi| J(|\xi|)|(D_{e}^{{\epsilon_k}|\xi|}u^{\epsilon_k}\cdot e)^-|^2\,d\xi\,dx
=\int_{D\times\mathcal{H}_1(0)\setminus K^{+,\epsilon_k}}|\xi| J(|\xi|)|D_{e}^{{\epsilon_k}|\xi|}u^{\epsilon_k}\cdot e|^2\,d\xi\,dx\nonumber\\
&&\leq C\int_{D\times\mathcal{H}_1(0)\setminus K^{+,\epsilon_k}}|\xi| J(|\xi|)F_{\epsilon_k |\xi|}(D_{e}^{\epsilon_k |\xi|}u^{\epsilon_k}\cdot e)\,d\xi\,dx\leq C\sup_{t\in [0,T]}\sup_{\epsilon_k}PD^{\epsilon_k}(u^{\epsilon_k})
\label{third}
\end{eqnarray}
where Theorem \ref{Gronwall} implies that the right most element of  the sequence of inequalities is bounded and \eqref{L2bound} follows.

To establish \eqref{L1bound} we apply H\"older's inequality to find that
\begin{eqnarray}
&&\int_{D\times\mathcal{H}_1(0)}|\xi| J(|\xi|)|D_{e}^{{\epsilon_k}|\xi|}u^{\epsilon_k}\cdot e|\,d\xi\,dx\nonumber\\
&&=\int_{K^{+,\epsilon_k}}|\xi| J(|\xi|)|D_{e}^{{\epsilon_k}|\xi|}u^{\epsilon_k}\cdot e|\,d\xi\,dx+\int_{D\times\mathcal{H}_1(0)\setminus K^{+,\epsilon_k}}|\xi| J(|\xi|)|D_{e}^{{\epsilon_k}|\xi|}u^{\epsilon_k}\cdot e|\,d\xi\,dx\nonumber\\
&&\leq \frac{2\Vert u^{\epsilon_k}\Vert_{L^\infty(D;\mathbb{R}^d)}}{\epsilon_k}\int_{K^{+,\epsilon_k}}|\xi|J(|\xi|)\,d\xi\,dx+\nonumber\\
&&+\nu(D\times\mathcal{H}_1(0))^{\frac{1}{2}}\left (\int_{D\times\mathcal{H}_1(0)}|\xi| J(|\xi|)|(D_{e}^{{\epsilon_k}|\xi|}u^{\epsilon_k}\cdot e)^-|^2\,d\xi\,dx\right)^{\frac{1}{2}},
\label{twoterms}
\end{eqnarray}
and \eqref{L1bound} follows from \eqref{suppsetupper} and   \eqref{L2bound}.

We establish \eqref{Fprimesecond}. This bound follows from the basic features of the potential function $f$. We will recall for subsequent use that $f$ is smooth positive, concave and $f'$ is a decreasing function with respect to its argument. So for $A$ fixed and $0\leq h\leq A^2\overline{r}^2$ we have
\begin{eqnarray}
|f'(h)-f'(0)|\leq |f'(A^2\overline{r}^2)- f'(0)|<2|f'(0)|^2.
\label{ffact}
\end{eqnarray}
The bound \eqref{Fprimesecond} is now shown to be a consequence of the following upper bound valid for the parameter $0<A<1$ given by
\begin{eqnarray}
&&\int_{D\times\mathcal{H}_1(0)}|\xi|J(|\xi|)|f'(\epsilon_k |\xi||(D_{e}^{\epsilon_k |\xi|}u^{\epsilon_k}\cdot e)^-|^2)-f'(0)|^2\, d\xi\,dx\nonumber\\
&&\leq \nu(D\times\mathcal{H}_1(0))\times|f'(A^2\overline{r}^2)-f'(0)|^2+C\epsilon_k\frac{4|f'(0)|^2}{A^2}.
\label{usefulbound}
\end{eqnarray}
We postpone the proof of \eqref{usefulbound} until after it is used to establish  \eqref{Fprimesecond}. Set $h_{\epsilon_k}=(D_{e}^{\epsilon_k |\xi|}u^{\epsilon_k}\cdot e)^-$ to note 
\begin{eqnarray}
F_{\epsilon_k |\xi|}'(h_{\epsilon_k})-2f'(0)h_{\epsilon_k}=(f'(\epsilon_k |\xi| h^2_{\epsilon_k})-f'(0))2h_{\epsilon_k}.
\label{diffeq}
\end{eqnarray}
Applying H\"olders inequality, \eqref{Fprimefirst}, \eqref{L2bound}, \eqref{usefulbound}, and \eqref{diffeq} gives
\begin{eqnarray}
&&\left\vert\int_{D\times\mathcal{H}_1(0)}|\xi| J(|\xi|)F_{\epsilon_k|\xi|}'(D_{e}^{{\epsilon_k}|\xi|}u^{\epsilon_k}\cdot e) (D_{e}^{\epsilon_k}\phi\,d\xi\cdot e)\,dx
\right\vert\nonumber\\
&&\leq\left\vert\int_{D\times\mathcal{H}_1(0)}|\xi| J(|\xi|)F_{\epsilon_k|\xi|}'((D_{e}^{{\epsilon_k}|\xi|}u^{\epsilon_k}\cdot e)^+) (D_{e}^{\epsilon_k}\phi\cdot e)\,d\xi\,dx
\right\vert\nonumber\\
&&+\left\vert\int_{D\times\mathcal{H}_1(0)}|\xi| J(|\xi|)F_{\epsilon_k|\xi|}'((D_{e}^{{\epsilon_k}|\xi|}u^{\epsilon_k}\cdot e)^-) (D_{e}^{\epsilon_k}\phi\cdot e)\,d\xi\,dx
\right\vert\nonumber\\
&&\leq C\sqrt{\epsilon_k}\Vert\mathcal{E}\phi\Vert_{\scriptscriptstyle{L^\infty(D;\mathbb{R}^{d\times d})}}+2\int_{D\times \mathcal{H}_1(0)}|\xi|J(|\xi|)f'(0)(D_e^{\epsilon_k |\xi|}u^{\epsilon_k}\cdot e)^-(D_e^{\epsilon_k |\xi|}\phi\cdot e)\,d\xi\,dx\nonumber\\
&&+\int_{D\times\mathcal{H}_1(0)}|\xi|J(|\xi|)\left (F_{\epsilon_k|\xi|}'((D_{e}^{{\epsilon_k}|\xi|}u^{\epsilon_k}\cdot e)^-) -2f'(0)(D_e^{\epsilon_k |\xi|}u^{\epsilon_k}\cdot e)^-\right )(D_e^{\epsilon_k |\xi|}\phi\cdot e)\,d\xi\,dx\nonumber\\
&&\leq C\left(f'(0)+\sqrt{\epsilon_k}+\left(\nu(D\times\mathcal{H}_1(0))\times|f'(A^2\overline{r}^2)-f(0)|^2+\epsilon_k\frac{4|f'(0)|^2}{A^2}\right)^{1/2}\right )\Vert\mathcal{E} \phi\Vert_{\scriptscriptstyle{L^\infty(D;\mathbb{R}^{d\times d})}}.\nonumber\\
\label{five}
\end{eqnarray}
and \eqref{Fprimesecond} follows.

We establish the inequality \eqref{usefulbound}. Set $h_{\epsilon_k}=(D_{e}^{\epsilon_k |\xi|}u^{\epsilon_k}\cdot e)^-$ and for $0<A<1$ introduce the set
\begin{eqnarray}
K^{+,\epsilon_k}_A=\{(x,\xi)\in D\times\mathcal{H}_1(0)\,: A^2\overline{r}^2\leq \epsilon_k|\xi||h_{\epsilon_k}|^2\}.
\label{suppAset}
\end{eqnarray}
To summarize  $(x,\xi)\in K^{+,\epsilon_k}_A$ implies $A^2\overline{r}^2\leq\epsilon_k|\xi||h_{\epsilon_k}|^2\leq\overline{r}^2$ and $(x,\xi)\not\in K^{+,\epsilon_k}_A$ implies $\epsilon_k|\xi||h_{\epsilon_k}|^2<A^2\overline{r}^2$ and $|f'(\epsilon_k|\xi||h_{\epsilon_k}|^2)-f'(0)|\leq|f'(A^2\overline{r}^2)-f'(0)|$. Inequality \eqref{L2bound} implies
\begin{eqnarray}
&&C>\int_{K^{+,\epsilon_k}_A} |\xi|J(|\xi|) h_{\epsilon_k}^2\,d\xi\,dx\geq\frac{A^2\overline{r}^2}{\epsilon_k}\int_{K^{+,\epsilon_k}_A} J(|\xi|) \,d\xi\,dx\nonumber\\
&&\geq\frac{A^2\overline{r}^2}{\epsilon_k}\int_{K^{+,\epsilon_k}_A} |\xi|J(|\xi|) \,d\xi\,dx,
\label{chebyA}
\end{eqnarray}
the last inequality follows since $1\geq|\xi|>0$. Hence
\begin{eqnarray}
\int_{K^{+,\epsilon_k}_A} |\xi|J(|\xi|) \,d\xi\,dx\leq C\frac{\epsilon_k}{A^2\overline{r}^2},
\label{chebyAUpper}
\end{eqnarray}
and it follows that
\begin{eqnarray}
&&\int_{K^{+,\epsilon_k}_A} |\xi|J(|\xi|)|f'(\epsilon_k|\xi|h_{\epsilon_k}|^2-f'(0)|^2 \,d\xi\,dx\nonumber\\
&&\leq 4|f'(0)|^2\int_{K^{+,\epsilon_k}_A} |\xi|J(|\xi|) \,d\xi\,dx\leq C\epsilon_k\frac{4|f'(0)|^2}{A^2\overline{r}^2}.
\label{kepsplus}
\end{eqnarray}
Collecting observations gives
\begin{eqnarray}
&&\int_{D\times\mathcal{H}_1(0)\setminus K^{+,\epsilon_k}_A}|\xi|J(|\xi|)|f'(\epsilon_k |\xi||(D_{e}^{\epsilon_k |\xi|}u^{\epsilon_k}\cdot e)^-|^2)-f'(0)|^2\, d\xi\,dx\nonumber\\
&&\leq \nu(D\times\mathcal{H}_1(0))\times |f'(A^2\overline{r}^2)-f'(0)|^2,
\label{2ndbd}
\end{eqnarray}
and \eqref{usefulbound} follows.

We now prove Lemma \ref{twolimitsA}. Write
\begin{eqnarray}
F_{\epsilon_k|\xi|}'(D_{e}^{\epsilon_k |\xi|}u^{\epsilon_k}\cdot e)=F_{\epsilon_k|\xi|}'((D_{e}^{\epsilon_k |\xi|}u^{\epsilon_k}\cdot e)^+)
+F_{\epsilon_k|\xi|}'((D_{e}^{\epsilon_k |\xi|}u^{\epsilon_k}\cdot e)^-),
\label{plusminus}
\end{eqnarray}
and from \eqref{Fprimefirst} it follows that
\begin{eqnarray}
&&\lim_{\epsilon_k\rightarrow 0} \int_{D}\int_{\mathcal{H}_1(0)}|\xi|J(|\xi|)F_{\epsilon_k|\xi|}'(D_{e}^{{\epsilon_k}|\xi|}u^{\epsilon_k}\cdot e)D_{e}^{\epsilon_k|\xi|}\phi\cdot e\,d\xi\,dx\nonumber\\
&&=\lim_{\epsilon_k\rightarrow 0} \int_{D}\int_{\mathcal{H}_1(0)}|\xi|J(|\xi|)F_{\epsilon_k|\xi|}'((D_{e}^{{\epsilon_k}|\xi|}u^{\epsilon_k}\cdot e)^-)(D_{e}^{\epsilon_k|\xi|}\phi\cdot e)\,d\xi\,dx.
\label{functiotolimitminus}
\end{eqnarray}

To finish the proof we identify the limit of the right hand side of \eqref{functiotolimitminus}.
Set $h_{\epsilon_k}=(D_{e}^{\epsilon_k |\xi|}u^{\epsilon_k}\cdot e)^-$ and apply H\'older's inequality to find
\begin{eqnarray}
&&\int_{D\times\mathcal{H}_1(0)}|\xi|J(|\xi|)\left(F_{\epsilon_k|\xi|}'(h_{\epsilon_k}) -2f'(0)h_{\epsilon_k}\right)(D_e^{\epsilon_k|\xi|}\phi\cdot e)\,d\xi\,dx\nonumber\\
&&\leq C\int_{D\times\mathcal{H}_1(0)}|\xi|J(|\xi|)\left|F_{\epsilon_k|\xi|}'(h_{\epsilon_k}) -2f'(0)h_{\epsilon_k}\right|\,d\xi\,dx\Vert\mathcal{E}\phi\Vert_{\scriptscriptstyle{L^\infty(D;\mathbb{R}^{d\times d})}}
\label{firstestimate}
\end{eqnarray}
We estimate the first factor in \eqref{firstestimate} and
apply \eqref{diffeq}, H\"older's inequality, \eqref{L2bound}, and \eqref{usefulbound} to obtain
\begin{eqnarray}
&&\int_{D\times\mathcal{H}_1(0)}|\xi|J(|\xi|)\left |F_{\epsilon_k|\xi|}'(h_{\epsilon_k}) -2f'(0)h_{\epsilon_k}\right |\,d\xi\,dx\nonumber\\
&&\leq\int_{D\times\mathcal{H}_1(0)}|\xi|J(|\xi|)\left |f'(\epsilon_k |\xi||h_{\epsilon_k}|^2) -2f'(0)\right |\left | h_{\epsilon_k}\right |\,d\xi\,dx
\nonumber\\
&&\leq C\left(\nu(D\times\mathcal{H}_1(0))\times |f'(A^2\overline{r}^2)-f'(0)|^2+\epsilon_k\frac{4|f'(0)|^2}{A^2\overline{r}^2}\right)^{1/2}.
\label{usefulLemmaA}
\end{eqnarray}
Lemma \ref{twolimitsA} follows on applying the bound \eqref{usefulLemmaA}  to \eqref {firstestimate} and  passing to the $\epsilon_k$  zero limit and noting that the choice of $0<A<1$ is arbitrary.

We now prove Lemma \ref{twolimitsB}.
For $\tau>0$ sufficiently small define $K^\tau\subset D$  by  $K^{\tau}=\{x\in D:\,dist(x,J_{u^0(t)})<\tau\}$.  From Hypothesis \ref{remark555} the collection of centroids in $CPZ^\delta(\overline{r},1/2,1/2,t)$ lie inside $K^\tau$ for $\delta$ sufficiently small. (Otherwise the components of the collection $CPZ^\delta(\overline{r},1/2,1/2,t)$ would concentrate about a component of $CPZ^0(\overline{r},1/2,1/2,t)$ outside $K^\tau$; contradicting the hypothesis that $J_{u^0(t)}=CPZ^0(\overline{r},1/2,1/2,t)$). The collection of all points belonging to unstable neighborhoods associated with centroids in   is easily seen to be contained in the slightly larger set $K^{\tau,\delta}=\{x\in \,D; dist(x,K^\tau)<\delta\}$. From Hypothesis \ref{remark556} we may choose test functions $\phi\in C_0^1(D\setminus K^{\tau,2\delta})$ such that for $\epsilon_k$ sufficiently small
\begin{eqnarray}
(D_e^{\epsilon_k |\xi|}u^{\epsilon_k}\cdot e)^-\phi=(D_e^{\epsilon_k|\xi|}u^{\epsilon_k}\cdot e)\phi.
\label{identical}
\end{eqnarray}

We form the test functions $\phi(x)\psi(\xi)$, with $\phi\in C_0^1(D\setminus K^{\tau,2\delta})$ and 
$\psi\in C(\mathcal{H}_1(0))$. From \eqref{L2bound} we may pass to a subsequence  to find that $(D_e^{\epsilon_k |\xi|}u^{\epsilon_k}\cdot e)^-$ weakly converges to the limit 
$g(x,\xi)$ in $L^2(D\times\mathcal{H}_1(0);\nu)$. 
With this in mind we write
\begin{eqnarray}
&&\int_{D\times\mathcal{H}_1(0)} g(x,\xi)\phi(x)\psi(\xi)\,d\nu=\int_{D\times\mathcal{H}_1(0)} g(x,\xi)\phi(x)\psi(\xi)|\xi|J(|\xi|)\,d\xi\,dx\nonumber\\
&&=\lim_{\epsilon_k\rightarrow 0}\int_{D\times\mathcal{H}_1(0)}(D_e^{\epsilon_k |\xi|}u^{\epsilon_k}(x)\cdot e)^- \phi(x)\psi(\xi)|\xi|J(|\xi|)\,d\xi\,dx\nonumber\\
&&=\lim_{\epsilon_k\rightarrow 0}\int_{D\times\mathcal{H}_1(0)}(D_e^{\epsilon_k |\xi|}u^{\epsilon_k}(x)\cdot e) \phi(x)\psi(\xi)|\xi|J(|\xi|)\,d\xi\,dx\nonumber\\
&&=\lim_{\epsilon_k\rightarrow 0}\int_{D\times\mathcal{H}_1(0)}(u^{\epsilon_k}(x)\cdot e)(D_{-e}^{\epsilon_k |\xi|}\phi(x))\psi(\xi)|\xi|J(|\xi|)\,d\xi\,dx,\label{middlepart}
\end{eqnarray}
where we have integrated by parts using \eqref{intbypts}  in the last line of \eqref{middlepart}.
Noting that $D_{-e}^{\epsilon_k |\xi|}\phi(x)$ converges uniformly to $-e\cdot\nabla\phi(x)$ and from the strong convergence of $u^{\epsilon_k}$ to $u^0$ in $L^2$ we obtain
\begin{eqnarray}
&&=\lim_{\epsilon_k\rightarrow 0}\int_{D\times\mathcal{H}_1(0)}(u^{\epsilon_k}(x)\cdot e)(D_{-e}^{\epsilon_k |\xi|}\phi(x))\psi(\xi)|\xi|J(|\xi|)\,d\xi\,dx\nonumber\\
&&=-\int_{D\times\mathcal{H}_1(0)}(u^{0}(x)\cdot e) (e\cdot\nabla\phi(x))\psi(\xi)|\xi|J(|\xi|)\,d\xi\,dx\nonumber\\
&&=-\sum_{j,k=1}^d\int_{D}u^0_j(x)\,\partial_{x_k}\left(\phi(x) \int_{\mathcal{H}_1(0)} e_j e_k\psi(\xi)|\xi|J(|\xi|)\,d\xi\right)\,dx\nonumber\\
&&=\int_{D}\mathcal{E} u^0_{jk}(x)\left(\phi(x) \int_{\mathcal{H}_1(0)} e_j e_k\psi(\xi)|\xi|J(|\xi|)\,d\xi\right)\,dx\nonumber\\
&&=\int_{D\times{\mathcal{H}_1(0)} }(\mathcal{E} u^0(x) e\cdot e)\phi(x)\psi(\xi)|\xi|J(|\xi|)\,d\xi\,dx,
\label{identifyweakl2}
\end{eqnarray}
where we have made use of  $E u^0\lfloor D\setminus K^{\tau,2\delta}=\mathcal{E} u^0\,dx$ on the third line of \eqref{identifyweakl2}.
From the density of the span of the test functions we conclude that $g(x,\xi)=\mathcal{E} u^0(x)e\cdot e$ almost everywhere on $D\setminus K^{\tau,2\delta}\times\mathcal{H}_1(0)$.  Since $K^{\tau,2\delta}$ can be chosen to have arbitrarily small measure with vanishing $\tau$ and $\delta$ we conclude that  $g(x,\xi)=\mathcal{E} u^0(x) e\cdot e$ on $D\times\mathcal{H}_1(0)$ a.e. and Lemma \ref{twolimitsB} is proved.

\subsubsection{ Proof of Theorems \ref{epsiloncontropprocesszone} and \ref{bondunstable}}
\label{proofbondunstable}
We begin with the proof on the upper bound on the size of the process zone given by  Theorem \ref{epsiloncontropprocesszone}.
The set $K^{+,\epsilon_k}_\alpha$ is defined by \

\begin{eqnarray}
K^{+,\epsilon_k}_\alpha=\{(x,\xi)\in D\times\mathcal{H}_1(0);\,|D_e^{\epsilon_k|\xi|}u^{\epsilon_k}\cdot e|>\underline{k}|\epsilon_k\xi|^{\alpha-1}\}
\label{equivdescralpha}
\end{eqnarray}
where $0<\underline{k}\leq\overline{r}$ and $1/2\leq \alpha<1$. We set $\beta=2\alpha-1$ to see 
for $(x,\xi)\in K^{+,\epsilon_k}$  that 
\begin{eqnarray}
\epsilon_k|\xi||D_e^{\epsilon_k|\xi|}u^{\epsilon_k}(x)\cdot e|^2>\underline{k}^2\epsilon_k^\beta|\xi|^\beta.
\label{inequalbasic}
\end{eqnarray}
and recall that the potential function $f(r)=F_1(r)$ is increasing to get
\begin{eqnarray}
J(|\xi|)\frac{1}{\epsilon_k}F_{1}(\underline{k}^2\epsilon_k^\beta|\xi|^\beta)=|\xi|J(|\xi|)\frac{1}{\epsilon_k|\xi|}F_{1}(\underline{k}^2\epsilon_k^\beta|\xi|^\beta)
\leq|\xi|J(|\xi|)F_{\epsilon_k|\xi|}(D_{e}^{{\epsilon_k}|\xi|}u^{\epsilon_k}\cdot e)\label{firstalpha}
\end{eqnarray}
and in addition since $|\xi|^{1-\beta}\leq 1$ we have
\begin{eqnarray}
&& \int_{K^{+,\epsilon_k}_\alpha}\frac{1}{\epsilon_k}F_{1}(\underline{k}^2\epsilon_k^\beta|\xi|^\beta)|\xi|^{1-\beta}J(|\xi|)\,d\xi\,dx\leq \int_{K^{+,\epsilon_k}_\alpha}\frac{1}{\epsilon_k}F_{1}(\underline{k}^2\epsilon_k^\beta|\xi|^\beta)J(|\xi|)\,d\xi\,dx\nonumber\\
&&\leq\int_{K^{+,\epsilon_k}_\alpha} |\xi|J(|\xi|)F_{\epsilon_k|\xi|}(D_{e}^{{\epsilon_k}|\xi|}u^{\epsilon_k}\cdot e)\,d\xi\,dx\leq\sup_{t\in [0,T]}\sup_{\epsilon_k}PD^{\epsilon_k}(u^{\epsilon_k})
\label{firstbalpha}
\end{eqnarray}
Taylor approximation gives
\begin{eqnarray}
\frac{1}{\epsilon_k}F_1(\underline{k}^2\epsilon_k^\beta|\xi|^\beta)=\epsilon_k^{\beta-1}\underline{k}^2(f'(0)+o(\epsilon_k^\beta))|\xi|^\beta
\label{taylor}
\end{eqnarray}
Substitution of \eqref{taylor} into the left hand side of  \eqref{firstbalpha} gives
\begin{eqnarray}
\int_{K^{+,\epsilon_k}_\alpha} |\xi|J(|\xi|)\,d\xi\,dx\leq 
\frac{\epsilon_k^{1-\beta}}{\underline{k}^2(f'(0)+o(\epsilon_k^\beta))}\left(\sup_{t\in [0,T]}\sup_{\epsilon_k}PD^{\epsilon_k}(u^{\epsilon_k})\right).
\label{semifinal}
\end{eqnarray}
Introduce the characteristic function $\chi_\alpha^{+,\epsilon_k}(x,\xi)$ defined on $D\times\mathcal{H}_1(0)$ taking the value $1$ for $(x,\xi)$ in $K_\alpha^{+,\epsilon_k}$ and zero otherwise.
Observe that 
\begin{eqnarray}
\frac{1}{m}\int_{\mathcal{H}_1(0)}\chi_\alpha^{+,\epsilon_k}(x,\xi)|\xi|J(|\xi|)\,d\xi=P(\{y\in\mathcal{H}_{\epsilon_k}(x);|S^{\epsilon_k}(y,x)|>\underline{k}|y-x|^{\alpha-1}\}),
\label{equivlence}
\end{eqnarray}
so
\begin{eqnarray}
&&\int_DP(\{y\in\mathcal{H}_{\epsilon_k}(x);|S^{\epsilon_k}(y,x)|>\underline{k}|y-x|^{\alpha-1}\})dx\nonumber\\
&&\leq \frac{\epsilon_k^{1-\beta}}{m\times \underline{k}^2(f'(0)+o(\epsilon_k^\beta))}\left(\sup_{t\in [0,T]}\sup_{\epsilon_k}PD^{\epsilon_k}(u^{\epsilon_k})\right).
\label{semifinaler}
\end{eqnarray}
For $0<\overline{\theta}\leq 1$,  Tchebyshev's inequality delivers
\begin{eqnarray}
&& \mathcal{L}^d\left(\{x\in D;\, P\left(\{y\in\mathcal{H}_{\epsilon_k}(x);\,|S^{\epsilon_k}(y,x)|>\underline{k}|y-x|^{\alpha-1}\}\right)>\overline{\theta}\}\right)\nonumber\\
&&\leq \frac{\epsilon_k^{1-\beta}}{m\times \overline{\theta}\underline{k}^2(f'(0)+o(\epsilon_k^\beta))}\left(\sup_{t\in [0,T]}\sup_{\epsilon_k}PD^{\epsilon_k}(u^{\epsilon_k})\right).
\label{probnewvariableboundchebyshcv}
\end{eqnarray}
and Theorem \ref{epsiloncontropprocesszone} follows on applying \eqref{gineq} and \eqref{upperbound}.

Choose $\epsilon_k=\frac{1}{2^k}$ and Theorem \ref{epsiloncontropprocesszone} imples $\mathcal{L}^d(PZ^{\epsilon_k}(\underline{k},\alpha,\overline{\theta},t))<C(\frac{1}{{2}^k})^{1-\beta}$, where $C$ is independent of $t\in [0,T]$.
The collection of process zones for $\epsilon_k<\delta$ is written as
\begin{eqnarray}
CPZ^{\delta}(\underline{k},\alpha,\overline{\theta},t)=\cup_{\epsilon_k<\delta}PZ^{\epsilon_k}(\underline{k},\alpha,\overline{\theta},t)
\label{unstablecollection}
\end{eqnarray}
and from the geometric series we find
\begin{eqnarray}
\mathcal{L}^d\left(CPZ^{\delta}(\underline{k},\alpha,\overline{\theta},t)\right)<C{\delta}^{1-\beta}.
\label{boundcdelta}
\end{eqnarray}
Theorem \ref{bondunstable} follows noting further that $CPZ^{0}(\underline{k},\alpha,\overline{\theta},t)\subset CPZ^{\delta}(\underline{k},\alpha,\overline{\theta},t)\subset CPZ^{\delta'}(\underline{k},\alpha,\overline{\theta},t)$, for $0<\delta<\delta'$.

\setcounter{equation}{0} \setcounter{theorem}{0} \setcounter{lemma}{0}\setcounter{proposition}{0}\setcounter{remark}{0}\setcounter{remark}{0}
\setcounter{definition}{0}\setcounter{hypothesis}{0}

\section{Dynamics and limits of energies that $\Gamma$-converge to Griffith fracture energies}
\label{Ninth}
In this final section we collect ideas and illustrate how the approach presented in the earlier sections can be used to examine limits of dynamics associated with other energies that $\Gamma$- converge to the Griffith fracture energy. As an example we consider the phase field aproach based on the Ambrosio-Tortorelli approximation for  dynamic brittle fracture calculations  \cite{BourdinLarsenRichardson}. This model is seen to be a well posed formulation in the sense that existence of solutions can be shown \cite{LarsenOrtnerSuli}. 
To  illustrate the ideas we focus on anti-plane shear and 
the model is described by an out of plane  elastic displacement $u^\epsilon(x,t)$ and phase field $0\leq v^\epsilon(x,t)\leq 1$ defined for points $x$ belonging to the domain $D\subset \mathbb{R}^2$. The potential energy associated with the cracking body is  given by the Ambrosio-Tortorelli potential 
\begin{eqnarray}
P^\epsilon(u^\epsilon(t),v^\epsilon(t))=E^\epsilon(u^\epsilon(t),v^\epsilon(t))+H^\epsilon(v^\epsilon(t)),
\label{ab}
\end{eqnarray}
with
\begin{eqnarray}
E^\epsilon(u^\epsilon(t),v^\epsilon(t))=\frac{\mu}{2}\int_D a^\epsilon(t)|\nabla u^\epsilon(t)|^2\,dx
\label{detailsforab}
\end{eqnarray}
and
\begin{eqnarray}
H^\epsilon(v^\epsilon(t))=\frac{\mathcal{G}_c}{2}\int_D \frac{(1-v^\epsilon(t))^2}{2\epsilon}+2\epsilon|\nabla v^\epsilon(t)|^2\,dx.
\label{moredetailsforab}
\end{eqnarray}
here $a^\epsilon(t)=a^\epsilon(x,t)=(v^\epsilon(x,t))^2+\eta^\epsilon$, with $0<\eta^\epsilon<<\epsilon$. In this model the phase field $v^\epsilon$, provides an approximate description of a freely propagating crack taking the value $1$ for points $(x,t)$ away from the crack and zero on the crack.
To formulate the problem we introduce the space $H^1_0(D)$ defined to be displacements $u$ in $H^1(D)$ with zero Dirichlet data on $\partial D$ and the set of functions  $H_1^1(D)$ defined to be functions $v$ in $H^1(D)$ for which $v=1$ on $\partial D$.
The total energy is given by
\begin{eqnarray}
\mathcal{F}(t;u^\epsilon,\partial_t u^\epsilon,v^\epsilon)=\frac{1}{2}\int_D|\partial_t u^\epsilon|^2\,dx+P^\epsilon(u^\epsilon,v^\epsilon)-\int_D f(t)u^\epsilon\,dx.
\label{eergytotala}
\end{eqnarray}
The body force $f(x,t)$ is prescribed and the displacement - phase field pair $(u^\epsilon,v^\epsilon)$  is a solution of the initial boundary value problem  \cite{LarsenOrtnerSuli} given by:
\begin{eqnarray}
\partial_{tt}^2 u^\epsilon-{\rm{div}}\left(a^\epsilon(t)\nabla(u^\epsilon-\partial_t u^\epsilon)\right)=f(t),&& \hbox{  in $D$},\nonumber\\
u^\epsilon=0 \hbox{  and  } v^\epsilon=1, &&\hbox{  on $\partial D$},\label{ibvp}
\end{eqnarray}
for $t\in (0,T]$ with initial conditions $u^\epsilon(0)=u^\epsilon_0$, $\partial_t u^\epsilon(0)=u^\epsilon_1\in H^1_0(D)$,  satisfying the crack stability condition
\begin{eqnarray}
P^\epsilon(u^\epsilon(t),v^\epsilon(t))\leq\inf\left\{P^\epsilon(u^\epsilon(t),v):\, v\in H^1_1(D),\,\,v\leq v^\epsilon(t)\right\}
\label{stability}
\end{eqnarray}
and energy balance
\begin{eqnarray}
\mathcal{F}(t;u^\epsilon,\partial_t u^\epsilon,v^\epsilon)=\mathcal{F}(0;u^\epsilon_0,u^\epsilon_1,v^\epsilon_0)-\int_0^t\int_D\,a^\epsilon|\nabla\partial_\tau u^\epsilon|^2\,d\tau-\int_0^t\int_D\,\partial_\tau f u^\epsilon\,dx\,d\tau,
\label{eergytotalbalance}
\end{eqnarray}
for every time $0\leq t\leq T$. Finally the  initial condition for the phase field is chosen such that $v^\epsilon(0)=v^\epsilon_0$ where $v^\epsilon_0\in H^1_1(D)$ satisfies the unilateral minimality condition \eqref{stability}. In this formulation the pair $u^\epsilon(t)$, $v^\epsilon(t)$ provides a regularized model for free crack propagation. Here the phase field tracks the evolution of the crack with $v^\epsilon=1$ away from the crack and $v^\epsilon=0$ in the crack set. 

For a body force $f(x,t)$ in $C^1([0,T];L^2(D))$ it is shown in \cite{LarsenOrtnerSuli} that there exists at least one trajectory
$(u^\epsilon,v^\epsilon)\in H^2((0,T);L^2(D))\cap W^{1,\infty}((0,T); H^1_0(D))\times W^{1,\infty}((0,T); H_1^1(D))$ satisfying \ref{ibvp} in the weak sense, i.e.
\begin{eqnarray}
\int_D\partial_{tt}^2 u^\epsilon\varphi\,dx+\int_D(a^\epsilon(t)\nabla(u^\epsilon-\partial_t u^\epsilon)\cdot\nabla\varphi\,dx=\int_D f(t)u^\epsilon\,dx,
\label{ibvpweak}
\end{eqnarray}
for all $\varphi$ in $H^1_0(D)$ for almost every $t$ in $(0,T]$, with $u^\epsilon(0)=u^\epsilon_0$, $\partial_t u^\epsilon(0)=u^\epsilon_1$, $v^\epsilon(0)=v^\epsilon_0$,
and such that \eqref{stability} and \eqref{eergytotalbalance} are satisfied for all times $0<t\leq T$. We have formulated the problem in a simplified setting to illustrate the ideas and note that this type of evolution is shown to exist for evolutions with more general boundary conditions and for displacements in two and three dimensions, see \cite{LarsenOrtnerSuli}. For future reference we call the pair $(u^\epsilon(t),v^\epsilon(t))$ a phase field fracture evolution.

In what follows we pass to the $\epsilon\rightarrow 0$ limit in the phase field evolutions to show existence of a limiting evolution with bounded linear elastic fracture energy. Here the limit evolution $u^0(t)$ is shown to take values in the space of special functions of bounded variation $SBV(D)$.  This space is well known and can be thought of as a specialization of the space SBD introduced earlier in Section \ref{sec5}  that is appropriate for the scalar problem treated here. For a treatment of SBV and its relation to fracture the reader is referred to \cite{AmbrosioBrades}.

Applying the techniques developed in previous sections it is possible to state and prove following theorem on the dynamics associated with the phase field evolutions  $(u^\epsilon,v^\epsilon)$ in the limit as $\epsilon\rightarrow 0$.
\begin{theorem} {\bf Sharp interface limit of phase field fracture evolutions.}\\
\label{epszerolimitofLOS}
Suppose for every $\epsilon>0$ that: (a) the potential energy of the initial data $(u_0^\epsilon,v_0^\epsilon)$ is uniformly bounded, 
ie. $\sup_{\epsilon>0}\left\{P^\epsilon(u^\epsilon_0,v_0^\epsilon)\right\}<\infty$, and that  (b) $\Vert u^\epsilon(t)\Vert_{L^\infty(D)}<C$ for all $0<\epsilon$ and $0\leq t\leq T$. Then on passing to a subsequence if necessary in the phase field fracture evolutions $(u^\epsilon,v^\epsilon)$  there exists an anti-plane displacement field $u^0(x,t)$ in $SBV(D)$ for all $t\in [0,T]$ such that $u^0\in C([0,T];L^2(D))$ and 
\begin{eqnarray}
\label{limsbv}
\lim_{\epsilon\rightarrow 0}\max_{0\leq t\leq T}\{\Vert u^\epsilon(t)-u^0(t)\Vert_{L^2(D)}\}=0
\end{eqnarray}
with
\begin{eqnarray}
\label{lefm}
GF(u^0)=\frac{\mu}{2}\int_D|\nabla u^0(t)|^2+{G_c}\mathcal{H}^1(J_{u^0(t)})<C 
\end{eqnarray}
for $0\leq t\leq T$.
\end{theorem}

For anti-plane shear deformations the energy $GF$ is a special form of the energy $LEFM$ introduced in Section \ref{sec5}.

The strategy we will use for proving  Theorem \ref{epszerolimitofLOS} is the same as the one developed in the proofs of Theorems \ref{LimitflowThm} and \ref{LEFMMThm}. This strategy can be made systematic and applied to evolutions associated with potential energies that $\Gamma$- converge to the Griffith fracture energy. It consists of four component parts: 
\begin{enumerate}
\item Constructing upper bounds on the kinetic and potential energy of the evolutions that hold uniformly for $0\leq t\leq T$ and $0<\epsilon$.
\item Showing compactness of the evolution $u^\epsilon(t)$ in $L^2(D)$ for each time  $0\leq t\leq T$.
\item Showing limit points of the sequence $\{u^\epsilon(t)\}$ belong to $SBD(D)$ (or  $SBV(D)$ as appropriate) for each time $0\leq t\leq T$.
\item $\Gamma$-convergence of the potential energies  to the Griffith energy $LEFM$ (or $GF$ as appropriate).
\end{enumerate}

Assume first that Parts 1 through 4 hold for the the phase field fracture evolution with potential energies $P^\epsilon$ given by \eqref{ab}. These are used as follows to prove Theorem \ref{epszerolimitofLOS}.
Part 1 is  applied as in \eqref{lip} to show that the sequence of evolutions $u^\epsilon(t)$  is uniformly Lipschitz continuous in time with respect to the $L^2(D)$ norm, i.e.
\begin{eqnarray}
\label{lipschitzat}
\Vert u^\epsilon(t_1)-u^\epsilon(t_2)\Vert_{L^2(D)}\leq K |t_1-t_2|
\end{eqnarray}
for $K$ independent of $\epsilon$ and for any $0\leq t_1<t_2\leq T$.
Part 2 together with \eqref{lipschitzat} and the Ascoli theorem imply the existence of a subsequence and a limit $u^0(x,t)\in C([0,T];L^2(D))$ such that the convergence \eqref{limsbv} holds.
Part 3 
shows that $u^0(x,t)$ belongs to $SBV(D)$ for every time in $[0,T]$. Part 4 together with Part 1 and the lower bound property of $\Gamma$-convergence shows that \eqref{lefm} holds and Theorem \ref{epszerolimitofLOS} follows.

We now establish Parts 1 through 4 for the dynamic phase field fracture evolution introduced in \cite{LarsenOrtnerSuli}. To obtain a uniform bound on the kinetic and potential energy differentiate both sides of the energy balance \eqref{eergytotalbalance} with respect to time to get
\begin{eqnarray}
\label{gronagain}
&&\frac{d}{dt}\left(\frac{1}{2}\int_D|\partial_t u^\epsilon(t)|^2\,dx+P^\epsilon(u^\epsilon(t),v^\epsilon(t))\right)-\frac{d}{dt}\int_D\,f(t) u^\epsilon\,dx\\
&&=-\int_D a^\epsilon|\nabla \partial_t u^\epsilon|^2\,dx-\int_D\partial_t f u^\epsilon\,dx.\nonumber
\end{eqnarray}
Manipulation and application of the identity $f\partial_t u^\epsilon=\partial_t(f u^\epsilon)-\partial_t f u^\epsilon$ to \eqref{gronagain} delivers the inequality
\begin{eqnarray}
\label{gronagainagain}
&&\frac{d}{dt}\left(\frac{1}{2}\int_D|\partial_t u^\epsilon(t)|^2\,dx+P^\epsilon(u^\epsilon(t),v^\epsilon(t))\right)\\
&&\leq \int_D f \partial_t u^\epsilon\,dx.\nonumber
\end{eqnarray}
Now set
\begin{eqnarray}
W^\epsilon(t)=\left(\frac{1}{2}\int_D|\partial_t u^\epsilon(t)|^2\,dx+P^\epsilon(u^\epsilon(t),v^\epsilon(t))\right)+1
\label{westat}
\end{eqnarray}
and proceed as in Section \ref{GKP} to get
\begin{eqnarray}
\left(\frac{1}{2}\int_D|\partial_t u^\epsilon(t)|^2\,dx+P^\epsilon(u^\epsilon(t),v^\epsilon(t))\right)  \leq \left(\int_0^t\Vert f(\tau)\Vert_{L^2(D)}\,d\tau +\sqrt{W^\epsilon(0)}\right )^2-1.
\label{gineqat}
\end{eqnarray}
Part 1 easily follows from \eqref{gineqat} noting that $\sup_{\epsilon>0}\{W^{\epsilon}(0)\}<\infty$ is a consequence of hypothesis (a) of Theorem \ref{epszerolimitofLOS}. For this example Parts 2 and 3 follow from the uniform bound of Part 1, hypothesis (b) of Theorem \ref{epszerolimitofLOS} and the well known compactness result for the Ambrosio Tortorelli functional, see for example the remark following Theorem 2.3 of \cite{Giacomini}. Part 4 is given by the Ambrosio-Tortorelli result \cite{AT} as expressed in Theorem 2.3 of \cite{Giacomini}.

\section{Conclusions}

The cohesive model for dynamic brittle fracture evolution  presented in this paper does not require extra constitutive laws such as a kinetic relation between crack driving force and crack velocity or a crack branching condition. Instead the evolution of the process zone together with the the fracture set is governed by one equation consistent with NewtonÕs second law given by \eqref{eqofmotion}. This is a characteristic feature of peridynamic models  \cite{Silling1}, \cite{States}. This evolution is intrinsic to the formulation and encoded into the nonlocal cohesive constitutive law. Crack nucleation criteria although not necessary to the formulation follow from the dynamics and are recovered here by viewing nucleation as a dynamic instability, this is similar in spirit to \cite{SillingWecknerBobaru} and the work of \cite{BhattacharyaDyal} for phase transitions. Theorem \ref{epsiloncontropprocesszone} explicitly shows how the size of the process zone is controlled by the radius of the horizon. This analysis shows that {\em the horizon size $\epsilon$ for cohesive dynamics is a modeling parameter} that can be calibrated according to the size of the process zone obtained from experimental measurements. The process zone is seen to concentrate on a set of zero volume  as the length scale of non-locality characterized by the radius of the horizon $\epsilon$ goes to zero, see Theorem \ref{bondunstable}. In this limit the dynamics is shown (under suitable hypotheses) to coincide with the simultaneous evolution of a displacement crack set pair. Here displacement satisfies the elastic wave equation for points in space-time away from the crack set. The shear and Lam\'e moduli together with the energy release rate are described in terms of moments of the nonlocal potentials.

\section{Acknowlegements}
\label{Acknowlegements}
The author would like to thank Stewart Silling, Richard Lehoucq and Florin Bobaru for stimulating and fruitful discussions.
This research is supported by NSF grant DMS-1211066, AFOSR grant FA9550-05-0008, and NSF EPSCOR Cooperative Agreement No. EPS-1003897 with additional support from the Louisiana Board of Regents.

\end{document}